\pdfoutput=1
\newcommand{\ExecuteBeforeCleveref}{%
  \usepackage[noend]{algpseudocode}%
} 
\documentclass{shinyart}

\usepackage[utf8]{inputenc}

\usepackage{shinybib}
\usepackage{enumitem} 
\setlist[enumerate]{label={(\roman*)}}
\usepackage{autonum}

\newcommand{\eps}{\varepsilon}
\renewcommand{\phi}{\varphi}
\newcommand{\N}{\mathbb{N}}
\newcommand{\R}{\mathbb{R}}
\newcommand{\Rbar}{\overline{\mathbb{R}}}
\newcommand{\calF}{\mathcal{F}}

\newcommand{\calL}{\mathcal{L}}
\newcommand{\calM}{\mathcal{M}}
\newcommand{\calP}{\mathcal{P}}
\newcommand{\calT}{\mathcal{T}}
\newcommand{\TV}{\mathrm{TV}}
\newcommand{\dual}[1]{\langle #1 \rangle}
\newcommand{\scalprod}[1]{\left( #1 \right)_{L^2(\Omega)}}
\newcommand{\norm}[1]{\| #1 \|}
\newcommand{\set}[2]{\left\{#1:#2\right\}}
\newcommand{\wkto}{\rightharpoonup}

\DeclareMathOperator{\dom}{\mathrm{dom}}

\DeclareMathOperator{\Id}{\mathrm{Id}}

\renewcommand{\div}{\operatorname{\mathrm{div}}}
\newcommand{\prox}{\mathrm{prox}}
\newcommand{\proj}{\mathrm{proj}}
\newcommand{\umin}{u_{\min}}
\newcommand{\umax}{u_{\max}}
\newcommand{\Uext}{\ensuremath{\hat U}}

\usepackage{graphicx}
\usepackage{pgfplots}
\pgfplotsset{compat=newest}
\pgfplotsset{plot coordinates/math parser=false}
\SetKw{KwTerminate}{terminate}
\SetKw{And}{and}
\SetKwFunction{pred}{Predictor}

\title{Total variation regularization of multi-material topology optimization}
\author{Christian Clason\thanks{Faculty of Mathematics, University Duisburg-Essen, 45117 Essen, Germany (\email{christian.clason@uni-due.de})}
    \and Florian Kruse\thanks{Institute of Mathematics and Scientific Computing, University of Graz, Heinrichstrasse 36, 8010 Graz, Austria
    (\email{florian.kruse@uni-graz.at})}
    \and Karl Kunisch\thanks{Institute of Mathematics and Scientific Computing, University of Graz, Heinrichstrasse 36, 8010 Graz, Austria, and Radon Institute, Austrian Academy of Sciences, Linz, Austria
(\email{karl.kunisch@uni-graz.at}).}}
\date{November 17, 2017}

\hypersetup{
    pdftitle={Total variation regularization of multi-material topology optimization},
    pdfauthor={Christian Clason, Florian Kruse, Karl Kunisch}
    pdfkeywords={topology optimization, total variation, convex analysis, non-smooth optimization, semi-smooth Newton method}
}

\begin{document}

\maketitle

\begin{abstract}
    This work is concerned with the determination of the diffusion coefficient from distributed data of the state. This problem is related to homogenization theory on the one hand and to regularization theory on the other hand. An approach is proposed which involves total variation regularization combined with a suitably chosen cost functional that promotes the diffusion coefficient assuming prespecified values at each point of the domain. The main difficulty lies in the delicate functional-analytic structure of the resulting nondifferentiable optimization problem with pointwise constraints for functions of bounded variation, which makes the derivation of useful pointwise optimality conditions challenging. To cope with this difficulty, a novel reparametrization technique is introduced. Numerical examples using a regularized semismooth Newton method illustrate the structure of the obtained diffusion coefficient.
\end{abstract}

\section{Introduction}
\label{sec:introduction}

In this paper we revisit a challenging problem in the calculus of variations given by
\begin{equation}\label{eq:intro_prob}\tag{PI}
    \left\{
        \begin{aligned}
            \min_{u\in \mathcal{U}} \frac{1}{2}\norm{y-z}^2_{L^2(\Omega)}&+ \mathcal{R}(u)\\
            \text{s.t.}\quad -\div (u\nabla y)&= f \text{ in } \Omega,\\
            y&=0 \text{ on } \partial \Omega,
        \end{aligned}
    \right.
\end{equation}
where $\mathcal{U}$ denotes the set of admissible controls and $\mathcal{R}$ stands for a regularization term. This problem represents the optimization-theoretic formulation of the problem of determining the optimal distribution $u$ of material in the domain $\Omega$ from data $z$. If the data are only available in distributed part $\omega \subsetneq \Omega$ of the domain, then the cost functional in \eqref{eq:intro_prob} can readily be adapted. Problem \eqref{eq:intro_prob} arises as the regularization of a coefficient inverse problem; if the focus is on the situation that $u(x)$ is supposed to assume only preferred values $u_i$ specific to different materials, it can also be considered as a topology optimization problem.

In the calculus of variation literature, different forms of \eqref{eq:intro_prob} have received a tremendous amount of attention.
For the particular choice that $\mathcal{R}$ is not present and
\begin{equation}\label{eq:def_U}
    {\mathcal{U}}= \{u \in L^\infty(\Omega):0 < \umin \le u(x) \le \umax \}
\end{equation}
for constants $\umin$ and $\umax$, it was shown in \cite{Murat:1977} that the problem may fail to have a solution. Historically, this goes along with the development of homogenization theory and deep analytical concepts such as H-convergence and compensated compactness; see, e.g., \cite{Murat1997,Tartar:1987,Tartar:2009}. Such concepts allow associating a solution to \eqref{eq:intro_prob} without the use of a regularization term $\mathcal{R}$.

Here we follow a different perspective and aim for a formulation that allows numerical realization; in such a context the use of regularization terms provides a powerful tool. The goal must be to choose a functional $\mathcal{R}$ that guarantees existence to \eqref{eq:intro_prob} and at the same time does not affect the sought parameter $u$ too much. The use of a regularization term involving semi-norms of Sobolev spaces would conflict with this second requirement, since such a choice would prevent jumps of $u$ across hypersurfaces -- a property that we want to retain here. 

The choice for ${\mathcal{R}}$ that we propose and investigate in this paper is
\begin{equation}
    \mathcal{R}(u) = \alpha G(u) + \beta \TV(u),
\end{equation}
where $G$ is a pointwise ``multi-bang'' penalty as in \cite{CK:2013,CK:2015} that promotes the attainment of the predefined states $\{u_i\}_{i=1}^m$ almost everywhere, and $\TV$ denotes the total variation semi-norm. The use of $\TV$ will guarantee existence, while $G$ models the desired structural properties. The usefulness of $\TV$ has been established in the calculus of variations and in image analysis for several decades now; see, e.g., \cite{Ambrosio,Giusti,Attouch} and \cite{Rudin:1992,Chavent:1997}. It has also been used in topology optimization in \cite{Amstutz:2012} and \cite{Bourdin:2003}, but the approaches in these contributions are different from our formulation and do not contain the multi-material concept (although the latter considers a three-phase formulation with two different non-material phases, ``void'' and ``liquid''). Rather, this concept is an extension of our work from \cite{CK:2015}, where related topology optimization problems are considered in situations where well-posedness can be guaranteed without the need of employing $\TV$-regularization. Concerning approaches for multi-material topology optimization, we refer to, e.g., \cite{Amstutz2006,Amstutz2010,Amstutz2011,Blank:2014,Haslinger:2010}; among these, our ``multi-bang approach'' is most closely related to the second.
Finally, coefficient inverse problems have been studied in a wide variety of contexts.

The use of the $\TV$ functional entails an essential difficulty from an infinite dimensional optimization point of view. In fact, well-posedness of the PDE constraint in \eqref{eq:intro_prob} requires a strictly positive lower bound on $u$ as in the definition \eqref{eq:def_U} of $\mathcal{U}$. In the process of deriving optimality conditions, however, one is confronted with the problem of considering the subdifferential of $\TV(u) + I_{\mathcal{U}}$, where $I_{\mathcal{U}}$ denotes the indicator function of the set $\mathcal{U}$, e.g., as extended real-valued functions on $L^2(\Omega)$. In this case, the sum rule cannot be used to compute this subdifferential since neither of the two functionals $\TV$ and $I_{\mathcal{U}}$ is continuous at any point of its domain (which would be required to use a result as in \cite{Brezis:1986} on the sum of subdifferentials of convex functions). The fact that the sum rule is not applicable constitutes a major obstacle for deriving useful optimality conditions. Thus, we propose a different approach to ensure the well-posedness of the PDE constraint in \eqref{eq:intro_prob}:
We introduce a reparametrization of the coefficient in the PDE constraint which allows us to drop the explicit pointwise bounds in the definition of $\mathcal{U}$. This novel approach could be of interest also in situations different from the one considered in this work.

For the numerical solution, we consider a finite element discretization of the problem that allows deriving optimality conditions in terms of the expansion coefficients that, after introducing a Moreau--Yosida regularization of the multi-bang and total variation penalties, can be solved by a semismooth Newton-type method with path-following.

The paper is organized as follows. \Cref{sec:problem} contains the
problem statement, useful results on the state equation, and descriptions
of the transformation announced above, as well as of the multi-bang
penalty term. \Cref{sec:existence,sec:optsys}
are devoted to the existence of minimizers and first-order optimality conditions, respectively.
The discretization of the
infinite dimensional problem as well as of the optimality conditions
are provided in \cref{sec:numerical}. There we also provide a description of the
semismooth Newton-type method, employing dual regularizations of the
multi-bang penalty term and the $\TV$ term, which are needed for defining the
Newton steps. Numerical examples are provided in \cref{sec:examples} for two model problems motivated by the interpretation of \eqref{eq:intro_prob} as a topology optimization and a parameter identification problem, respectively.
Finally, in \cref{sec:appendix} we prove that strongly Lipschitz domains are regular in the sense of Gröger, an elementary but not completely obvious result that is important in our analysis. 

\section{Problem statement and preliminary results}\label{sec:problem}

We consider for $\alpha,\beta>0$ the following problem:
\begin{equation}\label{eq:problem}
    \left\{
        \begin{aligned}
            \min_{u\in BV(\Omega)} \frac{1}{2}\norm{y-z}^2_{L^2(\Omega)}+ &\alpha\, G(u) + \beta\, \TV(u) \\
            \text{s.t.}\quad -\div (\Phi_\eps(u)\nabla y)&= f \text{ in } \Omega,\\
            y&=0 \text{ on } \partial \Omega.
        \end{aligned}
    \right.
    \tag{P}
\end{equation}
Here, $\Omega\subset\R^d$, $d\in\N$, is a bounded strongly Lipschitz domain (see \cref{def:Lipboundary_Alt} for a rigorous definition),
$BV(\Omega)$ denotes the space of functions of bounded variation, and
$f\in L^2(\Omega)$ and $z\in L^2(\Omega)$ are given. Furthermore, $\TV$ denotes the total variation, $G$ is a multi-bang penalty, and $\Phi_\eps$ for $\eps\geq 0$ is a superposition operator defined by a (smoothed) pointwise projection onto the set $[\umin,\umax]\subset (0,\infty)$, each of which will be described in detail in the following subsections.

\subsection{Functions of bounded variation}\label{sec:BVbasics}

We recall, e.g., from \cite{Ambrosio,Giusti,Ziemer} that the space $BV(\Omega)$ is given by those functions $v\in L^1(\Omega)$ for which the distributional derivative $Dv$ is a Radon measure, i.e.,
\begin{equation}
    BV(\Omega) = \set{v\in L^1(\Omega)}{ \norm{Dv}_{\calM(\Omega)} < \infty}.
\end{equation}
The \emph{total variation} of a function $v\in BV(\Omega)$ is then given by
\begin{equation}
    \TV(v) := \norm{Dv}_{\calM(\Omega)} = \int_\Omega \mathrm{d}|Dv|_2,
\end{equation}
i.e., the total variation in the sense of measure theory of the vector measure $Dv\in\calM(\Omega;\R^d)=C(\overline\Omega;\R^d)^*$. Here, $|\cdot|_2$ denotes the Euclidean norm on $\R^d$; we thus consider here the \emph{isotropic} total variation.
For $v\in L^1(\Omega)\setminus BV(\Omega)$, we set $\TV(v)=\infty$.

The space $BV(\Omega)$ is a Banach space if equipped with the norm
\begin{equation}
    \norm{v}_{BV(\Omega)} := \norm{v}_{L^1(\Omega)} + \TV(v),
\end{equation}
see, e.g., \cite[Thm.~10.1.1]{Attouch}. Moreover, the space $C^\infty(\overline\Omega)$ is dense in $BV(\Omega)$ with respect to \emph{strict convergence}, i.e., for any $v\in BV(\Omega)$ there exists a sequence $\{v_n\}_{n\in\N}\subset C^\infty(\overline\Omega)$ such that
\begin{enumerate}
    \item $v_n \to v$ in $L^1(\Omega)$ and
    \item $\TV(v_n) \to \TV(v)$,
\end{enumerate}
see, e.g., \cite[Thm.~10.1.2]{Attouch}. In fact, a slight modification of the proof (which is based on approximation via mollification) shows that for $v\in BV(\Omega)\cap L^p(\Omega)$ with $1<p<\infty$, the convergence $v_n\to v$ in (i) holds even strongly in $L^p$ (since the constructed mollified sequence converges in $L^p$ for any $1\leq p<\infty$; see, e.g., \cite[Prop.~2.2.4]{Attouch}).

It follows that $BV(\Omega)$ embeds into $L^r(\Omega)$ continuously for every $r\in [1,\frac{d}{d-1}]$ and compactly if $r < \frac{d}{d-1}$, see, e.g., \cite[Cor.~3.49 together with Prop.~3.21]{Ambrosio}. Note that this requires $\Omega$ to be a strongly Lipschitz domain.
In addition, the total variation is lower semi-continuous with respect to strong convergence in $L^1(\Omega)$, i.e., if $\{u_n\}_{n\in\N}\subset BV(\Omega)$ and $u_n\to u$ in $L^1(\Omega)$, we have that
\begin{equation}\label{eq:tv_lsc}
    \TV(u)\leq \liminf_{n\to\infty} \TV(u_n),
\end{equation}
see, e.g., \cite[Thm.~5.2.1]{Ziemer}. Note that this does not imply that $\TV(u)<\infty$ and hence that $u\in BV(\Omega)$ unless $\{\TV(u_n)\}_{n\in\N}$ has a bounded subsequence.
From \eqref{eq:tv_lsc}, we also deduce that the convex extended real-valued functional $\TV:L^p(\Omega)\rightarrow\R\cup\{\infty\}$ is weakly lower semi-continuous for any $p\in [1,\infty]$.

\subsection{Multibang penalty}\label{sec:problem:multi-bang}

Let $u_1<\dots < u_m$ be a given set of desired coefficient values. Here we assume that $u_1=0$ and $u_m = \umax-\umin$ such that for $u(x)\in [u_1,u_m]$, we have $u(x)+\umin \in [\umin,\umax]$.
The multi-bang penalty $G$ is then defined similar to \cite{CK:2015}, where we have to replace the box constraints $u(x)\in [u_1,u_m]$ by a linear growth to ensure that $G$ is finite on $L^r(\Omega)$, $r<\infty$. Specifically, we consider
\begin{equation}
    G:L^1(\Omega)\to\R,\qquad G(u) = \int_{\Omega} g(u(x))\,\mathrm{d}x,
\end{equation}
where $g:\R\to \R$ is given by
\begin{equation}\label{eq:multi-bang_pw}
    g(t) = \begin{cases}
        - u_m t & t \leq u_1,\\
        \frac12 \left((u_{i}+u_{i+1})t - u_iu_{i+1}\right) & t \in [u_i,u_{i+1}],\quad 1\leq i <m,\\
        u_m t - \frac12 u_m^2 & t \geq u_m.
    \end{cases}
\end{equation}
It can be verified easily that $g$ is continuous (note that $u_1=0$), convex, and linearly bounded from above and below, i.e.,
\begin{equation}
    \tfrac12 u_2 |t| \leq g(t) \leq u_m |t|\qquad\text{for all }t\in \R.
\end{equation}
\begin{remark} 
    The definition of $g$ implies that $g(t)>g(0)=0$ for all $t\neq 0$ and
    that $g(t)>g(u_m)$ for all $t>u_m=\umax-\umin$.
    For the results of this section as well as of Sections 3 and 4, we only require these properties of $g$ rather than the specific form of $g$.
    In particular, the results also hold for $t\mapsto |t|$, i.e., if $G$ is replaced by the $L^1$ norm.
\end{remark}

Since $g$ is finite (and hence proper), convex, and continuous, the corresponding integral operator $G:L^r(\Omega)\to \R$ is finite, convex, and continuous (and hence \emph{a fortiori} weakly lower semi-continuous) for any $r\in [1,\infty]$, see, e.g., \cite[Prop.~2.53]{Barbu}.
Also, the properties of $g$ imply the following properties of $G$:
\begin{enumerate}[label={(\textsc{g}\arabic*)}]
    \item $G(v) > G(0)=0$ for all $v\in L^1(\Omega)\setminus\{0\}$, \label{ass:g:bound}
    \item $\frac12 u_2 \norm{v}_{L^1(\Omega)}\leq G(v) \leq u_m \norm{v}_{L^1(\Omega)}$ for all $v\in L^1(\Omega)$.\label{ass:g:equiv}
\end{enumerate}

\begin{figure}
    \centering
    \begin{subfigure}[t]{0.32\linewidth}
        \centering
        \begin{tikzpicture}[baseline,style={font=\small}]
            \begin{axis}[%
                width=\linewidth,
                xmin=-0.5,
                xmax=2.5,
                ymax=3,
                xlabel style={at={(axis cs:2.5,0)}},
                xlabel={$v$},
                xtick=\empty,
                extra x tick style={grid=major},
                extra x ticks={0,1,2},
                extra x tick labels={$u_1$,$u_2$,$u_3$},
                axis y line=left,
                axis x line=bottom,
                ]
                \addplot[%
                    domain=-0.5:2.5,
                    color=DarkBlue,solid,line width=1.5pt,
                    ]
                    {
                        (x<=0)*(-2*x)+%
                        and(x>0,x<=1)*(0.5*((0+1)*x-0*1))+%
                        and(x>1,x<=2)*(0.5*((1+2)*x-1*2))+%
                        (x>2)*(2*x-0.5*2*2)
                    };
                \addplot[%
                    domain=-0.5:2.5,
                    color=gray,dashed,line width=1pt,
                    ]
                    { 0.5*1*abs(x) };
                \addplot[%
                    domain=-0.5:2.5,
                    color=gray,dashdotted,line width=1pt,
                    ]
                    { 2*abs(x) };
            \end{axis}
        \end{tikzpicture}
        \caption{$g$}\label{fig:multi-bang:g}
    \end{subfigure}
    \hfill
    \begin{subfigure}[t]{0.32\linewidth}
        \centering
        \begin{tikzpicture}[baseline,style={font=\small}]
            \begin{axis}[%
                width=\linewidth,
                xmin=-0.5,
                xmax=2.5,
                ymin=-2.1,
                ymax=2.1,
                xlabel style={at={(axis cs:2.5,-2.1)}},
                xlabel={$v$},
                xtick=\empty,
                extra x tick style={grid=major},
                extra x ticks={0,1,2},
                extra x tick labels={$u_1$,$u_2$,$u_3$},
                axis y line=left,
                axis x line=bottom,
                ]
                \addplot[%
                    samples=2000,
                    domain=-0.5:2.5,
                    color=DarkBlue,solid,line width=1.5pt,
                    ]
                    {
                        (x<0)*(-2)+%
                        and(x>=0,x<1)*(0.5*(0+1))+%
                        and(x>=1,x<2)*(0.5*(1+2))+%
                        (x>=2)*(2)
                    };
            \end{axis}
        \end{tikzpicture}
        \caption{$\partial g$}\label{fig:multi-bang:dg}
    \end{subfigure}
    \hfill
    \begin{subfigure}[t]{0.32\linewidth}
        \centering
        \begin{tikzpicture}[baseline,style={font=\small}]
            \begin{axis}[%
                width=\linewidth,
                xmin=-2.3,
                xmax=2.3,
                ymin=-0.5,
                ymax=2.5,
                xlabel style={at={(axis cs:2.3,-0.5)}},
                xlabel={$q$},
                xtick=\empty,
                extra x tick style={grid=major},
                extra x ticks={-2,0.5,1.5,2},
                extra x tick labels={$-2$,$0.5$,$1.5$,$2$},
                axis y line=left,
                axis x line=bottom,
                ]
                \addplot[%
                    samples=2000,
                    domain=-2:2.001,
                    color=DarkBlue,solid,line width=1.5pt,
                    ]
                    {
                        (x<=-2)*(-0.5)+%
                        and(x>-2,x<=0.5)*(0)+%
                        and(x>0.5,x<=1.5)*(1)+%
                        and(x>1.5,x<2)*(2)+%
                        (x>=2)*(2.5)
                    };
            \end{axis}
        \end{tikzpicture}
        \caption{$\partial g^*$}\label{fig:multi-bang:dgstar}
    \end{subfigure}
    \caption{Pointwise multi-bang integrand $g$, subdifferential $\partial g$, and conjugate differential $\partial g^*$ ($u_1=0$, $u_2=1$, $u_3=2$)}
    \label{fig:multi-bang}
\end{figure}
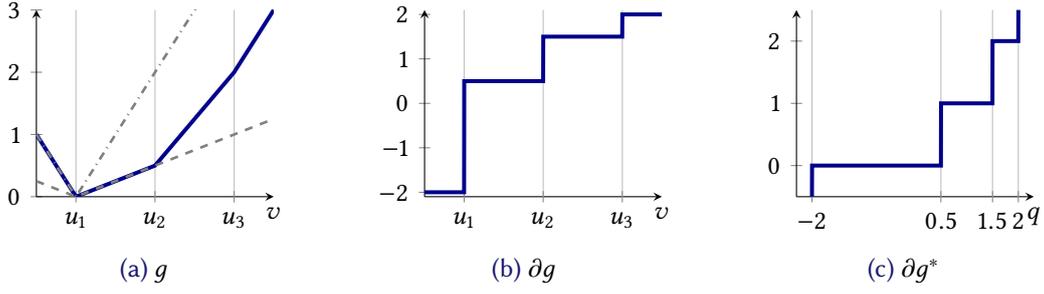

Furthermore, for $r<\infty$ and $r':=\frac{r}{r-1}$ (with $r'=\infty$ for $r=1$), the Fenchel conjugate
\begin{equation}
    G^*: L^{r'}(\Omega)\to \R\cup\{\infty\},\qquad G^*(q) = \sup_{v\in L^r(\Omega)} \dual{q,v}_{L^{r'}(\Omega),L^r(\Omega)}-G(v),
\end{equation}
as well as the convex subdifferential
\begin{equation}
    \partial G(v) = \set{q\in L^{r'}(\Omega)}{ \dual{q,\tilde v-v}_{L^{r'}(\Omega),L^r(\Omega)} \leq G(\tilde v)-G(v)\enspace\forall \tilde v\in L^r(\Omega)}
\end{equation}
can be computed pointwise, see, e.g., \cite[Props.~IV.1.2, IX.2.1]{Ekeland:1999a} and \cite[Prop.~2.53]{Barbu}, respectively. We point out that the pointwise representation of the subdifferential does \emph{not} hold for $r=\infty$.
From the definition of $g$ we thus obtain that
\begin{equation}\label{eq:g_subdiff}
    [\partial G(v)](x) \in \begin{cases}
        \{-u_m\} & v(x) < u_1,\\
        \left[-u_m,\tfrac12 (u_1+u_2)\right] & v(x) = u_1,\\
        \{\tfrac12(u_i+u_{i+1})\} & v(x)\in (u_i,u_{i+1}),\quad 1\leq i < m,\\
        \left[\tfrac12(u_{i-1}+u_{i}),\tfrac12(u_{i}+u_{i+1})\right] & v(x) = u_{i},\quad 1< i<m,\\
        \left[\tfrac12(u_{m-1}+u_m),u_m\right] & v(x) = u_m,\\
        \{u_m\} & v(x) > u_m,
    \end{cases}
\end{equation}
where, by a slight abuse of notation, $[\partial G(v)](x)$ stands for the evaluation of any $q\in\partial G(v)$ at $x\in\Omega$.
Using the fact that $s \in \partial g(t)$ if and only if $t\in \partial g^*(s)$ (see, e.g., \cite[Prop.~4.4.4]{Schirotzek:2007}), we deduce that
\begin{equation}\label{eq:gconj_subdiff}
    [\partial G^*(q)](x) \in \begin{cases}
        (-\infty,0] & q(x) = -u_m,\\
        \{0\} & q(x) \in \left(-u_m,\tfrac12(u_1+u_2)\right),\\
        [u_i,u_{i+1}] & q(x) = \tfrac12(u_{i}+u_{i+1}),\quad 1\leq i<m,\\
        \{u_i\} & q(x) \in \left(\tfrac12(u_{i-1}+u_{i}),\tfrac12(u_{i}+u_{i+1})\right), \quad 1< i<m,\\
        \{u_m\} & q(x) \in \left(\tfrac12(u_{m-1}+u_m),u_m\right),\\
        [u_m,\infty) & q(x) = u_m,\\
        \emptyset & \text{else},
    \end{cases}
\end{equation}
almost everywhere; see \cref{fig:multi-bang}.

\subsection{Superposition operator}

To ensure well-posedness of the state equation, both coercivity of the differential operator and pointwise boundedness of the coefficients are required. This can be achieved by imposing pointwise bounds on the coefficients. Appending such bounds to the problem statement \eqref{eq:problem} would lead to difficulties when deriving pointwise optimality conditions. As stated in the introduction, we therefore propose a reparametrization of the coefficient in the state equation. For this purpose we introduce the following family of (smoothed) pointwise projections onto the admissible set $[\umin,\umax]$.
For fixed $\eps \geq 0$ we consider $\phi_\eps:\R\to\R$,
\begin{equation}\label{eq:extsmo}
    \phi_\eps (t)=\umin +
    \begin{cases}
        -\eps& \text{for } t\leq -\eps,\\
        -\frac{1}{\eps^2}t^3-\frac{1}{\eps}t^2+t&\text{for }t\in[-\eps,0],\\
        t&\text{for } t\in [0,u_m],\\
        -\frac{1}{\eps^2}t^3 +\frac{3u_m+\eps}{\eps^2}t^2+\frac{\eps^2-2u_m\eps-3c^2}{\eps^2}t+\frac{u_m^3+c^2\eps}{\eps^2}&\text{for }t\in [u_m,u_m+\eps],\\
        u_m+\eps&\text{for } t\geq u_m+\eps,
    \end{cases}
\end{equation}
where we have used that $u_m=\umax-\umin$ from \cref{sec:problem:multi-bang}. For $\eps=0$, this coincides with the pointwise projection $\proj_{[\umin,\umax]}$, while for $\eps>0$ we have $\phi_{\eps} \in C^{1,1}(\R)$.
Clearly, there is a wide variety of choices which serves the purpose of making $\phi_\eps$ continuously differentiable. It is appropriate to choose this exterior smoothing in such a manner that $\phi_\eps'(t)\neq 0$ for $t\in [0,u_m]$. This will be further detailed in \cref{re:kk1} of \cref{sec:optsys}. The reader will notice in the following that $\eps>0$ is not used before deriving optimality conditions in \cref{sec:optsys}.

Since $\phi_\eps(t)$ is uniformly bounded and globally Lipschitz continuous, we deduce from \cite[Lem.~4.11]{Troeltzsch} that the corresponding superposition operator
\begin{equation}
    \Phi_\eps:L^r(\Omega) \to L^r(\Omega),\qquad [\Phi_\eps(v)](x) = \phi_\eps(v(x)) \quad \text{for a.e. }x\in\Omega,
\end{equation}
is globally Lipschitz continuous for every $r\in [1,\infty]$ and $\eps\geq 0$.

Similarly, for any $\eps>0$ it is easily verified that
\begin{equation}
    \phi_\eps' (t) =
    \begin{cases}
        -\frac{3}{\eps^2}t^2-\frac{2}{\eps}t+1&\text{for }t\in [-\eps,0]\\
        1&\text{for } t\in [0,u_m],\\
        -\frac{3}{\eps^2}t^2 +\frac{6u_m+2\eps}{\eps^2}t+\frac{\eps^2-2u_m\eps-3u_m^2}{\eps^2} &\text{for }t\in [u_m,u_m+\eps],\\
        0&\text{else},
    \end{cases}
\end{equation}
is locally Lipschitz continuous and uniformly bounded by $4/3$.
As a locally Lipschitz continuous function, $\phi_\eps'$ is even globally Lipschitz on the compact set $[-\eps,u_m+\eps]$. Since $\phi_\eps'(t)=0$ for all $t\in\R\setminus (-\eps,u_m+\eps)$, we infer that $\phi_\eps'$ is Lipschitz on all $\R$.
Hence, it follows from \cite[Lem.~4.12, proof of Lem.~4.13]{Troeltzsch} that $\Phi_\eps$ is
Lipschitz continuously Fréchet differentiable
from $L^\infty(\Omega)$ to $L^\infty(\Omega)$, and that the
Fréchet derivative
$\Phi_\eps'(v)\in \calL(L^\infty(\Omega),L^\infty(\Omega))$
at $v\in L^\infty(\Omega)$ acting on $h\in L^\infty(\Omega)$ is given by
\begin{equation}\label{eq:phi_derivative}
    [\Phi_\eps'(v)h](x) = \phi_\eps'(v(x))h(x) \quad \text{for a.e. }x\in\Omega.
\end{equation}
In particular, $\Phi_\eps'(v)$ can be represented pointwise almost everywhere by $x\mapsto \phi_\eps'(v(x))\in L^\infty(\Omega)$. In the following, we will not distinguish the derivative and its representation.
\begin{figure}
    \centering
    \begin{subfigure}[t]{0.45\linewidth}
        \centering
        \begin{tikzpicture}[baseline,style={font=\small}]
            \begin{axis}[%
                width=\linewidth,
                xmin=-0.5,
                xmax=2.5,
                ymin=-0.5,
                ymax=2.5,
                xlabel style={at={(axis cs:2.5,-0.1)}},
                xlabel={$t$},
                xtick=\empty,
                extra x tick style={grid=major},
                extra x ticks={-0.3,0,2,2.3},
                extra x tick labels={$-\eps$,$0$,$c$,$c+\eps$},
                axis y line=left,
                axis x line=bottom,
                ]
                \addplot[%
                    samples=100,
                    domain=-0.5:2.5,
                    color=DarkBlue,solid,line width=1.5pt,
                    ]
                    {
                        (x<=-0.3)*(-0.3)+%
                        and(x>-0.3,x<0)*(-1/(0.3*0.3)*x*x*x-1/0.3*x*x+x)+%
                        and(x>=0,x<=2)*(x)+%
                        and(x>2,x<2.3)*(-1/(0.3*0.3)*x*x*x+(3*2+0.3)/(0.3*0.3)*x*x+(0.3*0.3-2*2*0.3-3*2*2)/(0.3*0.3)*x+(2*2*2+2*2*0.3)/(0.3*0.3))+%
                        (x>=2.3)*(2.3)
                    };
            \end{axis}
        \end{tikzpicture}
        \caption{$\phi_\eps-u_{\min}$}\label{fig:superpos:phi}
    \end{subfigure}
    \hfill
    \begin{subfigure}[t]{0.45\linewidth}
        \centering
        \begin{tikzpicture}[baseline,style={font=\small}]
            \begin{axis}[%
                width=\linewidth,
                xmin=-0.5,
                xmax=2.5,
                ymin=-0.5,
                ymax=2.5,
                xlabel style={at={(axis cs:2.5,-0.1)}},
                xlabel={$t$},
                xtick=\empty,
                extra x tick style={grid=major},
                extra x ticks={-0.3,0,2,2.3},
                extra x tick labels={$-\eps$,$0$,$c$,$c+\eps$},
                axis y line=left,
                axis x line=bottom,
                ]
                \addplot[%
                    samples=100,
                    domain=-0.5:2.5,
                    color=DarkBlue,solid,line width=1.5pt,
                    ]
                    {
                        and(x>-0.3,x<0)*(-3/(0.3*0.3)*x*x-2/0.3*x+1)+%
                        and(x>=0,x<=2)*(1)+%
                        and(x>2,x<2.3)*(-3/(0.3*0.3)*x*x+2*(3*2+0.3)/(0.3*0.3)*x+(0.3*0.3-2*2*0.3-3*2*2)/(0.3*0.3))
                    };
            \end{axis}
        \end{tikzpicture}
        \caption{$\phi_\eps'$}\label{fig:superpos:phiprime}
    \end{subfigure}
    \caption{Smoothed projection $\phi_\eps$ and derivative $\phi_\eps'$ ($c:=u_m=2$, $\eps=0.3$)}
    \label{fig:superpos}
\end{figure}
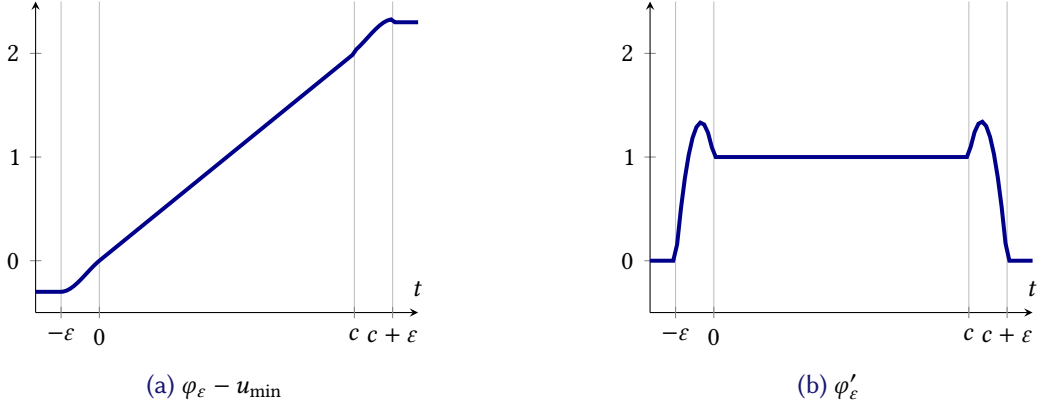

\subsection{State equation}

It will be convenient to introduce for $\eps\geq 0$ the set
\begin{equation}
    U_\eps=\set{v \in L^{\infty}(\Omega)}{0<\umin-\eps \leq v \leq \umax+\eps \text{ a.e. in } \Omega}
\end{equation}
along with its open $L^\infty(\Omega)$ neighborhood
\begin{equation}
    \Uext_\eps=\set{v\in L^{\infty}(\Omega)}{0<\tfrac12 \umin-2\eps < v < 2\umax+2\eps \text{ a.e. in } \Omega}.
\end{equation}
Furthermore, we consider for $w\in \Uext_\eps$ and $f\in L^2(\Omega)$ the elliptic partial differential equation
\begin{equation}\label{eq:state}
    \left\{
        \begin{aligned}
            -\div (w\nabla y)&= f \text{ in } \Omega,\\
            y&=0 \text{ on } \partial \Omega.
        \end{aligned}
    \right.
\end{equation}
From standard arguments based on the Lax--Milgram lemma, we obtain the existence of a unique solution $y\in H^1_0(\Omega)$ satisfying the uniform a priori estimate
\begin{equation}\label{eq:state:apriori_h1}
    \norm{y}_{H^1_0(\Omega)} \leq K_2 \norm{f}_{H^{-1}(\Omega)}
\end{equation}
for some $K_2>0$ independent of $w\in \Uext_\eps$ (but depending on $\Uext_\eps$), where $\norm{y}_{H^1_0(\Omega)}=\norm{\nabla y}_{L^2(\Omega)^d}$.
We also have the following global Lipschitz estimate for the solution mapping $w\mapsto y =:y(w)$.
\begin{lemma}\label{lem:state:lipschitz}
    For any $\eps\geq 0$ there exists a constant $L>0$ such that
    \begin{equation}
        \norm{y(w_1)-y(w_2)}_{H^1_0(\Omega)} \leq L\norm{w_1-w_2}_{L^\infty(\Omega)} \qquad \text{for all }w_1,w_2\in \Uext_\eps.
    \end{equation}
\end{lemma}
\begin{proof}
    Let $y_1,y_2\in H^1_0(\Omega)$ denote the solutions to \eqref{eq:state} for $w_1,w_2\in \Uext_\eps$, respectively. Inserting $y_1-y_2\in H^1_0(\Omega)$ as a test function in \eqref{eq:state} for $w=w_1$ and $w=w_2$, subtracting, inserting the productive zero, and rearranging yields
    \begin{equation}
        \left(w_1 \nabla(y_1-y_2),\nabla(y_1-y_2)\right)_{L^2(\Omega)^d}=
        \left((w_2-w_1) \nabla y_2,\nabla(y_1-y_2)\right)_{L^2(\Omega)^d}.
    \end{equation}
    Estimating the left-hand side using the uniform lower bound on $w$ and the right-hand side using the Cauchy--Schwarz inequality and the a priori estimate \eqref{eq:state:apriori_h1}, we obtain
    \begin{equation}
        \begin{aligned}
            (\tfrac12\umin-2\eps)\norm{\nabla (y_1-y_2)}^2_{L^2(\Omega)^d}
            &\leq \norm{w_1-w_2}_{L^\infty(\Omega)} \norm{\nabla y_2}_{L^2(\Omega)^d} \norm{\nabla (y_1-y_2)}_{L^2(\Omega)^d}\\
            &\leq K_2 \norm{f}_{H^{-1}(\Omega)}\norm{w_1-w_2}_{L^\infty(\Omega)} \norm{\nabla (y_1-y_2)}_{L^2(\Omega)^d},
        \end{aligned}
    \end{equation}
    from which the desired estimate follows with $L:=\frac{K_2}{\tfrac12\umin-2\eps}\norm{f}_{H^{-1}(\Omega)}$.
\end{proof}

Our next goal is to establish that there exists an $s>2$ such that the solution $y$ of \eqref{eq:state} belongs to $W^{1,s}(\Omega)$. This increase in regularity is crucial for obtaining pointwise optimality conditions.
The proof relies on results from Gröger \cite{Groeger}.
\begin{proposition}\label{thm:state:groeger}
    There exists an $s>2$ and a constant $K_s>0$ such that
    for all $w\in\Uext_\eps$ the solution $y\in H^1_0(\Omega)$ of \eqref{eq:state} satisfies
    \begin{equation}
        \norm{y}_{W^{1,s}(\Omega)} \leq K_s \norm{f}_{W^{-1,s}(\Omega)}.
    \end{equation}
\end{proposition}
\begin{proof}
    Fix $w\in\Uext_\eps$ and $f\in L^2(\Omega)$ and denote by $y\in H_0^1(\Omega)$ the solution to \eqref{eq:state}. By the Sobolev embedding theorem, there exists an $\bar s>2$ such that $L^2(\Omega)$ is continuously embedded in $W^{-1,s_1}(\Omega)$ for all $s_1\in(2,\bar s]$. Furthermore, by \cref{prop:regularity:LipschitzimpliesGroegerregularity} the domain $\Omega$ is regular in the sense of Gröger. Hence, \cite[Thm.~3]{Groeger} implies that $\Omega\in R_{s_2}$ for some $s_2>2$ and thus by \cite[Lem.~1]{Groeger} for $s:=\min\{s_1,s_2\}>2$ as well.
    We therefore obtain from \cite[Thm.~1]{Groeger} 
    for any $q\in W^{-1,s}(\Omega)$ that the unique solution $\hat y\in H_0^1(\Omega)$ of
    \begin{equation}\label{eq:PDEwithshift}
        \left\{
            \begin{aligned}
                -\div (w\nabla \hat y) + \hat y & = q &&\text{ in } \Omega,\\
                \hat y&= 0 &&\text{ on } \partial\Omega,
            \end{aligned}
        \right.
    \end{equation}
    satisfies $\norm{\hat y}_{W^{1,s}(\Omega)}\leq K\norm{q}_{W^{-1,s}(\Omega)}$, where $K$ denotes a constant that depends on $\Uext_\eps$ but not on $w$, $\hat y$, or $q$. For the choice $q=y+f$ this yields $\norm{\hat y}_{W^{1,s}(\Omega)}\leq K(C\norm{y}_{L^2(\Omega)}+\norm{f}_{W^{-1,s}(\Omega)})$,
    where $C$ denotes the constant of the continuous embedding $L^2(\Omega)\hookrightarrow W^{-1,s}(\Omega)$.
    Using the continuous embedding $H^1_0(\Omega)\hookrightarrow L^2(\Omega)$ with constant $\hat C$, the a priori estimate \eqref{eq:state:apriori_h1}, and
    the continuous embedding $W^{-1,s}(\Omega)\hookrightarrow H^{-1}(\Omega)$ with constant $\bar C$,
    we obtain $\norm{\hat y}_{W^{1,s}(\Omega)}\leq K(C\hat C K_2\bar C\norm{f}_{W^{-1,s}(\Omega)}+\norm{f}_{W^{-1,s}(\Omega)})$. Since for fixed $\Uext_\eps$ all appearing constants are independent of $w$, the claim follows by noting that the choice of $q$ implies that $y$ solves \eqref{eq:PDEwithshift}, hence $\hat y=y$.
\end{proof}

\section{Existence}\label{sec:existence}

To show existence of a solution to \eqref{eq:problem}, we make use of the solution mapping $w\mapsto y(w)$ to introduce the reduced functional
\begin{equation}
    J:BV(\Omega)\to \R,\qquad J(u) = \frac12\norm{y(\Phi_\eps(u))-z}_{L^2(\Omega)}^2 +\alpha\,G(u)+\beta\,\TV(u).
\end{equation}
\begin{proposition}\label{thm:existence:bv}
    For every $\eps\geq 0$ there exists a global minimizer $\bar u\in BV(\Omega)$ to \eqref{eq:problem}.
\end{proposition}
\begin{proof}
    Since $J$ is bounded from below due to \ref{ass:g:bound}, there exists a minimizing sequence $\{u_n\}_{n\in\N}\subset BV(\Omega)$. Furthermore, by \ref{ass:g:equiv}, we may assume without loss of generality that there exists a $C>0$ such that
    \begin{equation}
        C\left(\norm{u_n}_{L^1(\Omega)} + \TV(u_n)\right) \leq J(u_n)\leq J(0) \quad\text{for all }n\in\N,
    \end{equation}
    and hence that $\{u_n\}_{n\in\N}$ is bounded in $BV(\Omega)$.
    By the compact embedding of $BV(\Omega)$ into $L^1(\Omega)$ for any $d\in\N$, we can thus extract a subsequence, denoted by the same symbol, converging strongly in $L^1(\Omega)$ to some $\bar u\in L^1(\Omega)$.
    Lipschitz continuity of $\Phi_\eps$ from $L^1(\Omega)$ to $L^1(\Omega)$
    now implies that $\Phi_\eps(u_n)\to \Phi_\eps(\bar u)$ in $L^1(\Omega)$ as well. Furthermore, the corresponding sequence $\{y(\Phi_\eps(u_n))\}_{n\in\N}$ is uniformly bounded in $H^1_0(\Omega)$ due to \eqref{eq:state:apriori_h1}, and hence there exists a $\bar y\in H^1_0(\Omega)$ such that, after passing to a further subsequence if necessary, $y(\Phi_\eps(u_n))\wkto \bar y$ in $H^1_0(\Omega)$.
    Since $\{\Phi_\eps(u_n)\}_{n\in\N}$ is uniformly bounded in $L^\infty(\Omega)$ by construction, we have that $\Phi_\eps(u_n)\to\Phi_\eps(\bar u)$ strongly in $L^r(\Omega)$ for any $r\in [1,\infty)$ and, in particular, for $r=2$. We can thus pass to the limit in the distributional formulation of \eqref{eq:state},
    \begin{align}
        \scalprod{\Phi_\eps(u_n),\nabla y(\Phi_\eps(u_n))\cdot \nabla \psi} = \scalprod{f,\psi}\qquad\text{for all }\psi\in C^\infty_0(\Omega),\\
        \intertext{to obtain}
        \scalprod{\Phi_\eps(\bar u),\nabla \bar y\cdot \nabla \psi} = \scalprod{f,\psi}\qquad\text{for all }\psi\in C^\infty_0(\Omega).
    \end{align}
    By density, we obtain that $\bar y = y(\Phi_\eps(\bar u))$ and hence that
    $y(\Phi_\eps(u_n))\to y(\Phi_\eps(\bar u))$ strongly in $L^2(\Omega)$.
    Finally, lower semi-continuity of $G$ and $\TV$
    with respect to convergence in $L^1(\Omega)$
    and the strong convergence $y(\Phi_\eps(u_n))\to y(\Phi_\eps(\bar u))$ in $L^2(\Omega)$ imply that
    \begin{equation}
        J(\bar u) \leq \liminf_{n\to\infty} J(u_n) \leq J(u) \qquad\text{for all }u\in BV(\Omega)
    \end{equation}
    and thus that $\bar u\in BV(\Omega)$ is the desired minimizer.
\end{proof}

Due to the bilinear structure of the state equation the optimal control is not unique. Nonetheless, as a consequence of the reparametrization of the control by means of $\Phi_\eps$, any solution to \eqref{eq:problem} automatically satisfies pointwise control constraints.
\begin{proposition}\label{thm:existence:linfty}
    Let $\eps\geq 0$ and $\bar u\in BV(\Omega)$ be a local solution to \eqref{eq:problem}. Then, $\bar u+\umin\in U_\eps\subset L^\infty(\Omega)$.
\end{proposition}
\begin{proof}
    Let $\eps\geq 0$ and $\bar u\in BV(\Omega)$ with $\bar u+\umin\notin U_\eps$. We will show that $\bar u$ is not a local solution to \eqref{eq:problem}.
    We start by comparing $\bar u$ to $\hat u$ defined pointwise almost everywhere by
    \begin{equation}
        \hat u(x) = \begin{cases}
            -\eps & \bar u(x) < -\eps,\\
            \bar u(x) & \bar u(x) \in [-\eps,u_m+\eps],\\
            u_m+\eps & \bar u(x) > u_m+\eps.
        \end{cases}
    \end{equation}
    By definition of $\phi_\eps$, it follows that $\Phi_\eps(\hat u) = \Phi_\eps(\bar u)$ and thus that $y(\Phi_\eps(\hat u))=y(\Phi_\eps(\bar u))$.

    Furthermore, from Stampacchia's Lemma for BV functions \cite[Lem.~2.5]{Scherzer:2001} we obtain that $\TV(\hat u)\leq \TV(\bar u)$.
    Using the pointwise definition of $G$ together with the inequalities
    $g(t)>g(-\eps)>0$ for all $t<-\eps$ and $g(t)>g(u_m+\eps)$ for all $t>u_m+\eps$, we also deduce that
    $G(\hat u) < G(\bar u)$
    since $\bar u+\umin\notin U_\eps$. Thus, $J(\hat u) < J(\bar u)$.
    Similarly, we observe that $y(\Phi_\eps(u_t))=y(\Phi_\eps(\bar u))$ for all $t\in[0,1]$,
    where we have denoted $u_t:=(1-t)\hat u+t\bar u$.
    Using $\TV(\hat u)\leq \TV(\bar u)$ and 
    $G(\hat u) < G(\bar u)$
    together with the convexity of $\TV$ and $G$ yields that $\TV(u_t)\leq\TV(\bar u)$ and $G(u_t)<G(\bar u)$ for all $t\in[0,1)$. It follows that $J(u_t)<J(\bar u)$ for all $t\in [0,1)$ and hence that $\bar u$ is not a local solution to \eqref{eq:problem}.
\end{proof}
By \cref{thm:existence:linfty}, for any $\eps\geq 0$, each locally optimal control to problem \eqref{eq:problem} is therefore also a local solution of
\begin{equation}
    \min_{u\in BV(\Omega)\cap L^\infty(\Omega)} J(u),
\end{equation}
and, moreover, the set of globally optimal controls is the same for both problems. In particular, the solutions $\bar u$ to \eqref{eq:problem} for $\eps=0$ coincide with the solutions to
\begin{equation}\label{eq:starproblem}
    \left\{
        \begin{aligned}
            \min_{u\in BV(\Omega)} \frac{1}{2}\norm{y-z}^2_{L^2(\Omega)} + \alpha\, G(u) &\,+ \beta\, \TV(u) \\
            \text{s.t.}\qquad\qquad\quad u(x)+\umin &\in [\umin,\umax],\\
            \text{and }\quad
            -\div ((u+\umin)\nabla y)&= f \text{ in } \Omega,\\
            y&=0 \text{ on } \partial \Omega,
        \end{aligned}
    \right.
    \tag{P$^*$}
\end{equation}
which is a particular case of the motivating problem \eqref{eq:intro_prob}.

\begin{remark}
    The same cut-off argument as in the proof of \cref{thm:existence:linfty} can be applied to the minimizing sequence in the proof of Proposition \ref{thm:existence:bv} to construct a minimizing sequence that is bounded in $L^\infty(\Omega)$ and hence in $L^1(\Omega)$ even for $\alpha=0$. We thus also obtain the existence of a solution $\bar u$ to \eqref{eq:starproblem} with $\alpha=0$. The results in the following \cref{sec:optsys} remain valid in this case, and the optimality conditions derived therein simplify in an obvious manner.
\end{remark}
\bigskip

We close this section by briefly addressing the convergence of global solutions to \eqref{eq:problem} as $\eps \to 0^+$. For this purpose we consider a family $\{\bar u_\eps\}_{\eps>0}$ of solutions to \eqref{eq:problem}. From \cref{thm:existence:linfty} and the fact that $J(0)$ is independent of $\eps$, we deduce that
this family is bounded in $L^{\infty}(\Omega)\cap BV(\Omega)$ as $\eps\to 0^+$. Thus, there exists a sequence $\{\bar u_{\eps_k}\}_{k\in\N}$ converging strongly to some $\bar u$ in $L^r(\Omega)$ for every $r\in [2,\infty)$ with $\TV(\bar u)\le \liminf_{k\to \infty}\TV(\bar u_{\eps_k})<\infty$. With some modifications (in particular using that for every $u\in BV(\Omega)$ there holds $\Phi_{\eps_k}(u)\to\Phi_0(u)=\proj_{[\umin,\umax]}(u)$ strongly in $L^1(\Omega)$ for $k\to \infty$), the proof of \cref{thm:existence:bv} can now be used to verify that $\bar u$ is a global solution to \eqref{eq:problem} for $\eps=0$ and thus for \eqref{eq:starproblem}.

\section{Optimality conditions}\label{sec:optsys}

In this section, we derive \emph{pointwise} necessary optimality conditions for solutions to problem \eqref{eq:problem}.
Since we will require differentiability of the control-to-state operator $u\mapsto y(\Phi_\eps(u))$, we have to assume $\eps>0$ from here on. To keep the presentation simple, we will from now omit the dependence on $\eps$.
The derivation rests crucially on the following two nontrivial properties:
\begin{enumerate}
    \item By \cref{thm:existence:linfty}, we can work in the $L^\infty(\Omega)$ topology rather than in the $L^{\frac{d}{d-1}}(\Omega)$ topology induced by $BV(\Omega)$, which allows differentiability of the forward mapping.
    \item By \cref{thm:state:groeger}, the derivative of the forward mapping is actually in $L^r(\Omega)$ for some $r>1$, which will yield multipliers in $L^{r}(\Omega)$ instead of $L^{\infty}(\Omega)^*$.
\end{enumerate}

We begin by showing differentiability of the reduced tracking term
\begin{equation}\label{eq:tracking}
    F:\Uext\rightarrow\R,\qquad F(w) = \frac12\norm{y(w)-z}_{L^2(\Omega)}^2.
\end{equation}
This can be argued from differentiability of the forward mapping $w\mapsto y(w)$ in $L^\infty(\Omega)$ (see, e.g., \cite{Roesch}) together with the chain rule. However, it actually holds under the weaker requirement of Lipschitz continuity of the forward mapping shown in \cref{lem:state:lipschitz}. Since this argument may be of independent interest, we give a full proof here.

We first introduce for a given parameter $w\in\Uext\subset L^\infty(\Omega)$ and $y\in H^1_0(\Omega)$ the adjoint equation
\begin{equation}\label{eq:adjoint}
    \left\{
        \begin{aligned}
            -\div (w\nabla p)&= -(y-z) \text{ in } \Omega,\\
            p&=0 \text{ on } \partial \Omega.
        \end{aligned}
    \right.
\end{equation}
By the same arguments as for the state equation \eqref{eq:state} there exists a unique solution $p=p(w,y)\in H^1_0(\Omega)$, which depends continuously on $y$ and for which the additional regularity $p(w,y)\in W^{1,s}(\Omega)$ from \cref{thm:state:groeger} holds.
\begin{lemma}\label{lem:tracking_frechet}
    The mapping $F$ defined in \eqref{eq:tracking} is Lipschitz continuously Fréchet differentiable in every $w\in \Uext\subset L^\infty(\Omega)$. Furthermore, the Fréchet derivative of $F$ in $w\in \Uext$ is given by
    \begin{equation}
        F'(w) = \nabla y(w)\cdot \nabla p(w) \in L^{\frac{s}{2}}(\Omega)
    \end{equation}
    with $s>2$ from \cref{thm:state:groeger},
    where $y(w)\in H^1_0(\Omega)$ is the solution to \eqref{eq:state} and $p(w):=p(w,y(w))\in H^1_0(\Omega)$ is the corresponding solution to \eqref{eq:adjoint}.
\end{lemma}
\begin{proof}
    We first show directional differentiability in $\Uext\subset L^\infty(\Omega)$. Let $w\in\Uext$ and $h\in L^\infty(\Omega)$. Then there exists a $\rho_0>0$ sufficiently small such that $w+\rho h\in\Uext$ for all $\rho\in(0,\rho_0)$. Consequently, for all such $\rho$ there exists a solution $y(w+\rho h)\in H^1_0(\Omega)$ to \eqref{eq:state}. We now insert the productive zero $y(w)-y(w)$ in $F(w+\rho h)$ and expand the square to obtain
    \begin{equation}\label{eq:tracking_frechet1}
        \begin{aligned}[t]
            F(w+\rho h)-F(w) &= \frac1{2}\norm{(y(w+\rho h)-y(w))+(y(w)-z)}_{L^2(\Omega)}^2 - \frac1{2}\norm{y(w)-z}_{L^2(\Omega)}^2\\
                             &= \frac1{2}\norm{y(w+\rho h)-y(w)}_{L^2(\Omega)}^2 + \scalprod{y(w+\rho h)-y(w),y(w)-z}.
        \end{aligned}
    \end{equation}
    For the first term, we can use \cref{lem:state:lipschitz} to estimate
    \begin{equation}\label{eq:tracking_frechet2}
        \frac1{2}\norm{y(w+\rho h)-y(w)}_{L^2(\Omega)}^2 \leq \frac{L^2}{2}\rho^2 \norm{h}_{L^\infty(\Omega)}^2.
    \end{equation}
    For the second term, we introduce the adjoint state $p(w)$, integrate by parts, and use the state equation \eqref{eq:state} for $y=y(w)$ and $y=y(w+\rho h)$ to obtain
    \begin{equation}
        \begin{aligned}
            \scalprod{y(w+\rho h)-y(w),y(w)-z} &= \scalprod{y(w+\rho h)-y(w),\div(w \nabla p)}\\
                                               &=\scalprod{\div (w\nabla y(w+\rho h)),p} -
            \scalprod{\div (w\nabla y(w)),p}\\
            &= \scalprod{-f,p} - \scalprod{\div (\rho h\nabla y(w+\rho h)),p} - \scalprod{-f,p}\\
            &= \rho\scalprod{h\nabla y(w+\rho h),\nabla p}.
        \end{aligned}
    \end{equation}
    By \cref{lem:state:lipschitz} we have that $y(w+\rho h)\to y(w)$ in $H^1_0(\Omega)$ as $\rho\to 0^+$. Hence, dividing \eqref{eq:tracking_frechet1} by $\rho>0$ and passing to the limit implies in combination with \eqref{eq:tracking_frechet2} that
    \begin{equation}
        F'(w;h) := \lim_{\rho\to0^+}\frac1\rho(F(w+\rho h) - F(w)) = \dual{h,\nabla y\cdot\nabla p}_{L^\infty(\Omega),L^1(\Omega)}.
    \end{equation}
    Since the mapping $h\mapsto F'(w;h)$ is linear and bounded, $\nabla y\cdot \nabla p$ is the Gâteaux derivative of $F$ at $w\in\Uext$. Thus, $F$ is Gâteaux differentiable in $\Uext$.
    Due to \cref{lem:state:lipschitz} the mappings $w\mapsto y(w)$ and $w\mapsto p(w,y)$ are Lipschitz from $L^\infty(\Omega)$ to $H^1_0(\Omega)$ in $\Uext$.
    By using \eqref{eq:state:apriori_h1}, we infer that the mapping $y\mapsto p(w,y)$ is Lipschitz from $H^1_0(\Omega)$ to $H^1_0(\Omega)$ for any fixed $w\in\Uext$, with a Lipschitz constant independent of $w$.
    This shows that $w\mapsto p(w):=p(w,y(w))$ is Lipschitz continuous from $L^\infty(\Omega)$ to $H^1_0(\Omega)$ in $\Uext$.
    Hence, the mapping $w\mapsto\nabla y(w)\cdot\nabla p(w)$ is Lipschitz continuous from $L^\infty(\Omega)$ to $L^1(\Omega)$ in $\Uext$, and thus $F$ is in fact Fréchet differentiable in $\Uext$ with Lipschitz continuous derivative.
    The regularity $ \nabla y(w)\cdot \nabla p(w)\in L^{\frac{s}{2}}(\Omega)$ follows from \cref{thm:state:groeger}.
\end{proof}

Together with the Fréchet differentiability of $\Phi$ in $L^\infty(\Omega)$, this allows deriving abstract first-order necessary optimality conditions using classical tools from convex analysis. Here it is crucial that $G$ does not incorporate pointwise constraints and is finite on $L^p(\Omega)$ for $p=\frac{s}{s-2}>1$ instead of $p=1$ in order to apply the sum rule to its convex subdifferential (considered as a subset of $L^q(\Omega)$ with $q=\frac{s}{2}<\infty$), which requires the effective domain of $G$ to have non-empty interior.
\begin{theorem}\label{thm:optimality}
    Any local minimizer $\bar u\in BV(\Omega)$ to \eqref{eq:problem} satisfies
    \begin{equation}\label{eq:optimality}
        -F'(\Phi(\bar u))\Phi'(\bar u) \in \alpha\,\partial G(\bar u) + \beta\,\partial \TV(\bar u)\subset L^\frac{s}{2}(\Omega),
    \end{equation}
    where $G$ and $\TV$ are considered as extended real-valued convex functionals on $L^\frac{s}{s-2}(\Omega)$.
\end{theorem}
\begin{proof}
    Let $\bar u\in BV(\Omega)$ be a local minimizer to \eqref{eq:problem}. \cref{thm:existence:linfty} shows that $\bar u$ is also a local minimizer in $BV(\Omega)\cap L^\infty(\Omega)$. Thus, for all $u\in BV(\Omega)\cap L^\infty(\Omega)$ and $t>0$ sufficiently small, we have that
    \begin{equation}
        F(\Phi(\bar u)) + \alpha\, G(\bar u) + \beta\, \TV(\bar u) \leq
        F(\Phi(\bar u+t(u-\bar u))) + \alpha\, G(\bar u+t(u-\bar u)) + \beta\, \TV(\bar u+t(u-\bar u)).
    \end{equation}
    We now proceed as in the proof of \cite[Prop.~2.2]{CK:2015}, using the convexity of $G$ and $\TV$ to obtain after rearranging that
    \begin{equation}
        \frac1t\left(F(\Phi(\bar u+t(u-\bar u))) - F(\Phi(\bar u))\right) + \alpha\left(G(u)-G(\bar u)\right)+\beta\left(\TV(u)-\TV(\bar u)\right)\geq 0.
    \end{equation}
    By \cref{lem:tracking_frechet} and the chain rule, $F\circ \Phi$ is Fréchet differentiable at $\bar u\in L^\infty(\Omega)$, and the Fréchet derivative is given by
    \begin{equation}
        (F\circ\Phi)'(\bar u) = F'(\Phi(\bar u))\Phi'(\bar u)\in L^\infty(\Omega)^*.
    \end{equation}
    Since \cref{lem:tracking_frechet} further implies that $F'(\Phi(\bar u))\in L^{\frac{s}2}(\Omega)$, and since we have $\Phi'(\bar u)\in L^\infty(\Omega)$ from the representation \eqref{eq:phi_derivative}, we deduce that in fact $(F\circ\Phi)'(\bar u)\in L^{\frac{s}2}(\Omega)\subset L^1(\Omega)$.
    Hence, we can pass to the limit $t\to 0^+$ to obtain
    \begin{equation}
        \dual{F'(\Phi(\bar u))\Phi'(\bar u), u-\bar u}_{L^1(\Omega),L^\infty(\Omega)}
        + \alpha\left(G(u)-G(\bar u)\right)+\beta\left(\TV(u)-\TV(\bar u)\right)\geq 0
    \end{equation}
    for all $u\in BV(\Omega)\cap L^\infty(\Omega)$.

    By the density of $C^\infty(\overline\Omega)$ in $L^{\frac{s}{s-2}}(\Omega)\cap BV(\Omega)$ with respect to strict convergence, there exists for any $u\in L^{\frac{s}{s-2}}(\Omega)\cap BV(\Omega)$ a sequence $\{u_n\}_{n\in\N}\subset C^\infty(\overline\Omega)$ with $u_n\to u$ strongly in $L^{\frac{s}{s-2}}$. Hence, $G(u_n)\to G(u)$ by continuity of $G$, $\TV(u_n)\to \TV(u)$, and
    \begin{equation}
        \dual{F'(\Phi(\bar u))\Phi'(\bar u), u_n-\bar u}_{L^1(\Omega),L^\infty(\Omega)} \to \dual{F'(\Phi(\bar u))\Phi'(\bar u), u-\bar u}_{L^{\frac{s}2}(\Omega),L^{\frac{s}{s-2}}(\Omega)}.
    \end{equation}
    Taking $\TV(u) = \infty$ for $u\in L^{\frac{s}{s-2}}(\Omega)\setminus BV(\Omega)$, we deduce that
    \begin{equation}
        \dual{F'(\Phi(\bar u))\Phi'(\bar u), u-\bar u}_{L^{\frac{s}2}(\Omega),L^{\frac{s}{s-2}}(\Omega)}
        + \alpha\left(G(u)-G(\bar u)\right)+\beta\left(\TV(u)-\TV(\bar u)\right)\geq 0
    \end{equation}
    holds for all $u\in L^{\frac{s}{s-2}}(\Omega)$. But this implies by definition that
    \begin{equation}
        -F'(\Phi(\bar u))\Phi'(\bar u) \in \partial(\alpha\, G+\beta \,\TV)(\bar u)\subset L^{\frac{s}{2}}(\Omega),
    \end{equation}
    where the subdifferentials are understood as those of the canonical restriction
    to $L^{\frac{s}{s-2}}(\Omega)$.

    Finally, since $\dom \TV = BV(\Omega)\cap L^\frac{s}{s-2}(\Omega)\subset L^\frac{s}{s-2}(\Omega) = \dom G$ and $G$ is continuous on $L^\frac{2}{s-2}(\Omega)$, we can apply the sum rule for convex subdifferentials (see, e.g., \cite[Prop.~4.5.1]{Schirotzek:2007}) to obtain \eqref{eq:optimality}.
\end{proof}

Introducing explicit subgradients for the two subdifferentials, we obtain primal-dual optimality conditions.
\begin{corollary}\label{cor:optsys_pd}
    For any local minimizer $\bar u\in BV(\Omega)$ to \eqref{eq:problem}, there exist $\bar q \in L^\frac{s}2(\Omega)$ and $\bar \xi \in L^{\frac{s}2}(\Omega)$ satisfying
    \begin{equation}\label{eq:optsys_pd}
        \left\{
            \begin{aligned}
                0 &= F'(\Phi(\bar u))\Phi'(\bar u) + \alpha\bar q + \beta \bar \xi,\\
                \bar q &\in \partial G(\bar u),\\
                \bar \xi &\in \partial \TV(\bar u).
            \end{aligned}
        \right.
    \end{equation}
\end{corollary}

From \cref{cor:optsys_pd}, we can further derive \emph{pointwise} optimality conditions for optimal controls. For the Fréchet derivative of the tracking term and the subdifferential of the multi-bang penalty, we apply \cref{lem:tracking_frechet} together with the representations \eqref{eq:phi_derivative} and \eqref{eq:g_subdiff}, respectively. The characterization of $\bar \xi\in \partial \TV(\bar u)$ is more involved. Formally, elements of the subdifferential $\partial \TV(u)$ have the form
$-\div \left(\frac{\nabla u}{|\nabla u|_2}\right)$,
which is equal to the negative mean curvature of the level sets of $u$.
This can be made rigorous using the \emph{full trace} from \cite{BrediesHoller}, which requires some notation. First, we introduce for $1\leq q<\infty$ the space
\begin{equation}
    W^{\div,q}(\Omega):=\set{v\in L^q(\Omega;\R^d)}{\div v \in L^q(\Omega)}
\end{equation}
endowed with the graph norm. Furthermore, for any Radon measure $\mu$, let $L^1_\mu(\Omega;\R^d)$ denote the space of $\mu$-measurable functions $v:\Omega\to\R^d$ for which
\begin{equation}
    \norm{v}_{L^1_\mu(\Omega;\R^d)} := \int_\Omega |v(x)|_2 \,\mathrm{d}\mu
\end{equation}
is finite.
To any $v\in W^{\div,q}(\Omega)\cap L^\infty(\Omega)$, we can then assign a unique $Tv\in L^1_{|Du|}(\Omega;\R^d)$, called the \emph{full trace} of $v$, using appropriate converging sequences; see \cite[Def.~12]{BrediesHoller} for a precise definition. Finally, we recall the decomposition of the measure $Du$ for $u\in BV(\Omega)$ into an absolutely continuous part $D^au = \nabla u\, \mathrm{d}\mathcal{L}^d$ with respect to the $d$-dimensional Lebesgue measure $\mathcal{L}^d$, a jump part
\begin{equation}
    D^j u = (u^+ - u^-)\nu_u\, \mathrm{d}\mathcal{H}^{d-1}|_{S_u},
\end{equation}
where $u^+-u^-$ denotes the jump of $u$ on the singularity set $S_u$ with normal $\nu_u$ and $(d-1)$-dimensional Hausdorff measure $\mathcal{H}^{d-1}$, and the Cantor part $D^c u$ with density $\sigma_{u}$ with respect to $|D^cu|$. We can now state fully our pointwise optimality conditions.
\begin{theorem}\label{thm:opt_pointwise}
    For any local minimizer $\bar u\in BV(\Omega)$ to \eqref{eq:problem}, there exist $\bar y,\bar p \in W^{1,s}(\Omega)$, $\bar q \in L^\frac{s}2(\Omega)$, and $\bar \psi\in W^{\div,\frac{s}2}(\Omega)$ satisfying
    \allowdisplaybreaks
    \begin{subequations}\label{eq:kkt}
        \begin{align}
            &\left\{
            \begin{aligned}
                -\div (\Phi(\bar u)\nabla \bar y) &= f \quad\text{in }\Omega,\\
                \bar y &= 0\quad\text{on }\partial\Omega,
            \end{aligned}
        \right.
        \label{eq:kkt:state}\\[0.5em]
        &\left\{
        \begin{aligned}
            -\div (\Phi(\bar u)\nabla \bar p) &= -(\bar y-z) \quad\text{in }\Omega,\\
            \bar p &= 0\quad\text{on }\partial\Omega,
        \end{aligned}
    \right.
    \label{eq:kkt:adjoint}\\[0.5em]
    & (\nabla \bar y \cdot \nabla \bar p)\Phi'(\bar u) + \alpha \bar q - \beta \div \bar\psi = 0 \quad\text{in } L^{\frac{s}2}(\Omega),
    \label{eq:kkt:gradient} \\[0.5em]
    &\bar u(x) \in
    \begin{cases}
        (-\infty,u_1] & \bar q(x) = -u_m,\\
        \{u_1\} & \bar q(x) \in \left(-u_m,\tfrac12(u_1+u_2)\right),\\
        [u_i,u_{i+1}] & \bar q(x) = \tfrac12(u_{i}+u_{i+1}),\quad 1\leq i<m,\\
        \{u_i\} & \bar q(x) \in \left(\tfrac12(u_{i-1}+u_{i}),\tfrac12(u_{i}+u_{i+1})\right), \quad 1< i<m,\\
        \{u_m\} & \bar q(x) \in \left(\tfrac12(u_{m-1}+u_m),u_m\right),\\
        [u_m,\infty) & \bar q(x) = u_m,\\
        \emptyset & \text{else},
    \end{cases}
    \label{eq:kkt:multi-bang} \\[0.5em]
    &\left\{
    \begin{aligned}
        |\bar \psi(x)|_2 &\leq 1 \qquad&&\text {for a.e. } x\in \Omega,\\
        \bar \psi(x) &= \frac{\nabla \bar u(x)}{|\nabla \bar u(x)|_2}\qquad &&\text{for a.e. } x\in \Omega\text{ with } \nabla \bar u(x) \neq 0,\\
        (T\bar \psi)(x) &= \frac{\bar u^+(x)-\bar u^-(x)}{|\bar u^+(x)-\bar u^-(x)|}\nu_{\bar u}(x)\quad &&\text{for $\mathcal{H}^{d-1}$-a.e. } x\in S_{\bar u},\\
        (T\bar \psi)(x) &= \sigma_{\bar u}(x) \qquad &&\text{for $|D^c\bar u|$-a.e. } x\in \Omega.
    \end{aligned}
\right.
\label{eq:kkt:tv}
 \end{align}
\end{subequations}
\end{theorem}
\begin{proof}
    We start with \eqref{eq:kkt:gradient}, which is obtained from the first equation of \eqref{eq:optsys_pd} by using \cref{lem:tracking_frechet} to express $F'(\Phi(\bar u))\Phi'(\bar u)$ in terms of the solution $\bar y$ to the state equation \eqref{eq:kkt:state} and the solution $\bar p$ to the adjoint equation \eqref{eq:kkt:adjoint}.
    Furthermore, we have used \cite[Prop.~8]{BrediesHoller}, which states that any $\bar \xi\in\partial\TV(\bar u)\cap L^{q}(\Omega)$ can be expressed as $\bar \xi = -\div \bar\psi$ for a $\bar \psi\in W^{\div,q}(\Omega)$ satisfying \eqref{eq:kkt:tv}.\footnote{The result in \cite{BrediesHoller} is stated for $q=\frac{p}{p-1}$ for $1< p \leq \frac{d}{d-1}$. However, the upper bound on $p$ is not used in the proofs; it is merely the natural integrability of $u\in BV(\Omega)$ through embedding and is assumed to avoid further restrictions. We can thus apply the result for arbitrary $q>1$.}
    We point out that the $L^p(\Omega)$, $p>1$, regularity of $(F\circ \Phi)'(\bar u)$ is crucial to allow applying this result, and that it holds for strongly Lipschitz domains.
    Finally, the second relation in \eqref{eq:optsys_pd} can be equivalently written as $\bar u \in \partial G^*(\bar q)$, which by \eqref{eq:gconj_subdiff} admits the pointwise characterization \eqref{eq:kkt:multi-bang}.
\end{proof}

Let us briefly comment on these optimality conditions. Clearly, \eqref{eq:kkt:multi-bang} implies that if $\bar q$ does not have level sets of strictly positive measure, $\bar u$ will be a pure multi-bang control, i.e., $\bar u(x)\in \{u_1,\dots,u_m\}$ almost everywhere. Moreover, from \eqref{eq:kkt:tv} we can deduce that $\nabla\bar u(x)=0$ for almost every $x\in\Omega$ with $|\bar \psi(x)|_2<1$.
Further pointwise interpretations, in particular concerning the interaction between the multi-bang and the total variation penalty, is impeded by the fact that \eqref{eq:kkt:gradient} couples $\bar q$ not with $\bar \psi$ but with $\div\bar\psi$, and the divergence operator does not act pointwise and has a nontrivial kernel.

\begin{remark}\label{re:kk1}
    As already mentioned, the regularization $\phi_\eps$ of $\proj_{[\umin,\umax]}$ should be chosen in such a way that it does not become stationary in $[0,u_m]$. For example, if we define the function $\phi_\eps$ of \eqref{eq:extsmo} in such a manner that it acts as an interior smoothing with $\phi_\eps'(t)=0$ for $t \in (-\infty,0]\cup [u_m,\infty)$, then $\bar u \equiv 0$ with $\bar q \equiv 0$, $\bar{\psi}\equiv 0$ and $\bar y$, $\bar p$ computed from \eqref{eq:kkt:state} and \eqref{eq:kkt:adjoint} always provides a trivial solution to the optimality system. It could also be observed that this obstructs numerical algorithms.

    Similarly, $\phi_\eps'(u_m)=0$ would restrict in an undesired
    manner the possibility that $\phi_\eps(u(x))=\umax$. In fact, if $\bar u(x) = u_m$
    on a ball $B$ of radius $\rho > 0$, then $\alpha \bar q(x) = \beta \div
    \bar \psi(x)$ on $B$, where $\bar q(x) \in \left(\frac{1}{2}(u_{m-1}+ u_m),
    u_m\right]$ for almost every $x\in B$ and $|\psi(x)|_2 \le 1$ for almost every $x\in\Omega$. As a consequence, we have that
    \begin{equation}
        \frac{\alpha \pi^{\frac{d}{2}}\rho^d(u_{m-1}+u_m)}{2
        \Gamma(\frac{d}{2}+1)} < \alpha \int_B \bar q \, dx = \beta\int_B \div
        \bar \psi \, dx = \beta\int_{\partial B}\bar \psi \cdot n \, ds \le
        \frac{2 \beta \pi^{\frac{d}{2}}\rho^{d-1}}{\Gamma(\frac{d}{2})},
    \end{equation}
    where $n$ denotes the unit outer normal to $B$. Thus, $\bar u(x) = u_m$ cannot occur on sets that contain a ball $B$ of radius
    $\rho \ge \frac{4 \beta \Gamma(\frac{d}{2}+1)}{\alpha (u_{m-1} + u_m) \Gamma(\frac{d}{2})}
    = \frac{2\beta d}{\alpha (u_{m-1} + u_m)}$. Using the same argument for a general set $B$ to which the divergence theorem applies, we infer that $\bar u = u_m$ in $B$ necessitates
    $\frac{|B|}{|\partial B|}<\frac{2\beta}{\alpha (u_{m-1} + u_m)}$.
\end{remark}

\section{Numerical solution}\label{sec:numerical}

This section is concerned with the numerical computation of solutions to \eqref{eq:problem}. We proceed in several steps. First, we introduce in \cref{sec:discretization} a finite element discretization of \eqref{eq:problem}, for which we derive in \cref{sec:discrete_optimality} necessary optimality conditions in terms of the coefficients with respect to the finite element basis functions. These can be solved by a semismooth Newton-type method with path-following that is described in \cref{sec:semismooth}.

\subsection{Discretization}\label{sec:discretization}

We consider a finite element discretization of \eqref{eq:problem}. Let $\calT=\{\calT_h\}_{h>0}$ be a quasi-uniform triangulation of $\Omega$, which we assume in the following to be polyhedral for simplicity, consisting of triangular or tetrahedral elements $T$ with 
volume $|T|$. For later use, let us also introduce the notation $\calT_h=\{T_j\}_{j=1}^{N_{\calT_h}}$ for $h>0$, i.e., $\calT_h$ consists of $N_{\calT_h}$ elements that are denoted by $T_j$, $1\leq j\leq N_{\calT_h}$.

For the state and adjoint equation, we choose a conforming piecewise linear discretization, i.e., we set
\begin{equation}
    Y_h := \set{v_h\in C_0(\Omega)}{v_h|_T \in \calP_1\text{ for all }T\in\calT_h}.
\end{equation}
In $Y_h$ we use the standard nodal basis $\{\delta_i^{Y_h}\}_{i=1}^{N_{Y_h}}$ with respect to the vertices $x_i\in\R^d$, $1\leq i\leq N_{Y_h}$.
For any $v_h\in Y_h$, we denote by $\hat v_h\in\R^{N_{Y_h}}$ the coefficients of $v_h$ with respect to this basis. Defining $[v]_j$ to be the $j$-th component of a vector $v$, we can express this for $v_h\in Y_h$ as
$v_h = \sum_{i=1}^{N_{Y_h}} [\hat v_h]_i \delta_i^{Y_h}$.

The control is also discretized as continuous and piecewise linear, i.e., we set
\begin{equation}
    U_h := \set{u_h\in C(\overline\Omega)}{u_h|_T \in \calP_1\text{ for all }T\in\calT_h}.
\end{equation}
This choice -- as opposed to piecewise constants -- yields a convergent (nonconforming) discretization even for the isotropic total variation, see \cite{Casas1999,Bartels:2012}.
Again we use the standard nodal basis, denoted by $\{\delta_i^{U_h}\}_{i=1}^{N_{U_h}}$, and distinguish between $u_h\in U_h$ and its coefficient vector $\hat u_h\in\R^{N_{U_h}}$.

For $w_h\in U_h$, the discrete state equation reads
\begin{equation}
    \scalprod{w_h \nabla y_h,\nabla v_h} = \scalprod{f,v_h}\qquad\text{for all }v_h\in Y_h,
\end{equation}
and similarly for the discrete adjoint equation.
We denote the corresponding (symmetric) stiffness matrix by
$A_h(w_h)\in\R^{N_{Y_h}\times N_{Y_h}}$ and the mass matrix by
$M_h\in\R^{N_{Y_h}\times N_{Y_h}}$.

Since the discrete gradient of $u_h\in U_h$ should be piecewise constant, we introduce the space
\begin{equation}
    \Psi_h :=\set{\psi_h\in L^2(\Omega)^d}{\psi_h|_T \in \calP_0^d\text{ for all }T\in\calT_h}.
\end{equation}
In $\Psi_h$ we work with the basis of characteristic functions of $T\in\calT_h$, denoted by $\{\chi_i\}_{i=1}^{N_{\Psi_h}}$. For the coefficients of $\psi_h\in\Psi_h$ associated to $T\in\calT_h$, we write $[\hat\psi_h]_T \in\R^d$ and
assume that $\hat\psi_h\in\R^{N_{\Psi_h}}$ is ordered in the way
$[\hat\psi_h]_{T_j} = ([\hat\psi_h]_{(j-1)d+1},\ldots,[\hat\psi_h]_{jd})^T\in\R^{d}$ for $1\leq j\leq N_{\calT_h}$. This allows us to infer that $([\hat\psi_h]_{T_j})_{1\leq j\leq N_{\calT_h}} = \hat\psi_h$.
Moreover, let
$D_h\in\R^{N_{\Psi_h}\times N_{U_h}}$
denote the stiffness matrix arising from the bilinear form
\begin{equation}
    \scalprod{\nabla u_h,\psi_h} \qquad\text{for all }(u_h,\psi_h)\in U_h\times\Psi_h.
\end{equation}
We mention that $-D_h^T\in\R^{N_{U_h}\times N_{\Psi_h}}$ corresponds to the discrete divergence.
In the following, we assume that $D_h$ is ordered in the way $[D_h\hat u_h]_{T_j}= (D_{(j-1)d+1}\hat u_h,\ldots,D_{jd}\hat u_h)^T\in\R^{d}$, where $D_i$ denotes for $1\leq i\leq N_{\Psi_h}$ the $i$-th row of $D_h$. This allows us to infer that the Fréchet derivative of the mapping $\hat u_h\mapsto ([D_h\hat u_h]_{T_j})_j\in\R^{N_{\Psi_h}}$, $1\leq j\leq N_{\calT_h}$, is given by $D_h$.

The multi-bang penalty is approximated via mass lumping, i.e., we take
\begin{equation}
    G_h(\hat u_h) :=\sum_{i=1}^{N_{U_h}} d_i g([\hat u_h]_i),
\end{equation}
where $g:\R\to\R$ is given by \eqref{eq:multi-bang_pw} and $d_i := \int_\Omega \delta_i^{U_h}(x)\,d x$, see \cite{CHW:2012,Pieper,Trautmann}.
For later use, we also introduce the diagonal matrix $M^\ell_h\in\R^{N_{U_h}\times N_{U_h}}$ with entries $d_i$, which corresponds to a lumped mass matrix in $U_h$.
Similarly, the total variation is approximated by
\begin{equation}
    \TV_h(\hat u_h) := \sum_{T\in\calT_h} |[D_h \hat u_h]_T|_2.
\end{equation}
This is a correctly weighted discretization of the total variation since for all $u_h\in U_h$ there holds
\begin{equation}
    \TV(u_h)=\sum_{T\in\calT_h} |T| |\nabla u_h\lvert_T|_2=\sum_{i=1}^{N_{\Psi_h}} |\scalprod{\nabla u_h,\chi_i}|_2 = \sum_{T\in\calT_h} |[D_h\hat u_h]_T|_2 = \TV_h(\hat u_h).
\end{equation}
Note that by these definitions, $G_h$ and $TV_h$ are defined on $\R^{N_{U_h}}$, allowing us to apply convex analysis in the standard Euclidean topology.

The discrete problem now reads
\begin{equation}\label{eq:problem_h}
    \left\{\begin{aligned}
            \min_{\hat u_h\in\R^{N_{U_h}}} &\frac12\norm{y_h-z_h}_{L^2}^2 +
            \alpha\,G_h(\hat u_h) +\beta\,\TV_h(\hat u_h)\\
            \text{s.t.}\quad &A_h(\Phi_h(u_h))\hat y_h = M_h\hat f_h,
    \end{aligned}\right.
\end{equation}
where $z_h$ is the $L^2(\Omega)$ projection of $z$ onto $Y_h$ and thus $\frac12\norm{y_h-z_h}_{L^2}^2=\frac12 (\hat y_h-\hat z_h)^T M_h (\hat y_h-\hat z_h)$. Similarly, $f_h$ denotes the $L^2(\Omega)$ projection (or interpolation) of $f$ onto $Y_h$.
The existence of a solution
$\hat u_h^\ast\in\R^{N_{U_h}}$ to \eqref{eq:problem_h} then follows from standard arguments.

\subsection{Discrete optimality system and regularization}\label{sec:discrete_optimality}

We now derive numerically tractable optimality conditions for the discretized problem \eqref{eq:problem_h}, exploiting the fact that functional-analytic difficulties that had to be circumvented to obtain \eqref{eq:kkt} do not arise in the finite-dimensional setting. Specifically,
\begin{enumerate}[label=(\roman*)]
    \item we can consider $\eps=0$ or equivalently, by \cref{thm:existence:linfty}, the discrete analogue of \eqref{eq:starproblem}, thus eliminating the need for $\Phi_\eps$;
    \item as in \cite{CK:2013,CK:2015}, we can include the pointwise constraints in the definition of the multi-bang penalty $G$;
    \item applying the chain rule to the convex subdifferential of the discrete total variation directly yields an explicit componentwise relation.
\end{enumerate}
Hence, we replace \eqref{eq:problem_h} by
\begin{equation}\tag{P$_{\text{h}}$}\label{eq:problem_hhat}
    \left\{\begin{aligned}
            \min_{\hat u_h\in\R^{N_{U_h}}} &\frac12\norm{y_h-z_h}_{L^2}^2 +
            \alpha\,\hat G_h(\hat u_h) +\beta\,\TV_h(\hat u_h)\\
            \text{s.t.}\quad &A_h(u_h+\umin)\hat y_h = M_h\hat f_h
    \end{aligned}\right.
\end{equation}
for
\begin{equation}\label{def:Ghhat}
    \hat G_h(\hat u_h) :=\sum_{i=1}^{N_{U_h}} d_i \hat g([\hat u_h]_i),\qquad
    \hat g(t) = \begin{cases}
        \infty & t < u_1,\\
        \frac12 \left((u_{i}+u_{i+1})t - u_iu_{i+1}\right) & t \in [u_i,u_{i+1}],\quad 1\leq i <m,\\
        \infty & t > u_m.
    \end{cases}
\end{equation}
Proceeding as in the continuous case, we see that \eqref{eq:kkt:state} and \eqref{eq:kkt:adjoint} are replaced by their finite element approximation. Introducing for $y_h,p_h\in Y_h$ the vector
\begin{equation}
    \hat a_h(y_h,p_h) := \nabla y_h\cdot\nabla p_h \in \R^{N_{U_h}},
\end{equation}
we obtain analogously to \eqref{eq:optsys_pd} the primal-dual optimality conditions
\begin{equation}\label{eq:kkt:gradient_h}
    \left\{
        \begin{aligned}
            A_h(u_h^\ast+\umin)\hat y_h^\ast &= M_h\hat f_h,\\
            A_h(u_h^\ast+\umin)\hat p_h^\ast &= M_h(\hat z_h-\hat y_h^\ast),\\
            0&= \hat a_h(y_h^*,p_h^*) + \alpha \hat q_h^\ast + \beta \hat\xi_h^\ast,\\
            \hat q_h^\ast &\in \partial \hat G_h(\hat u _h^\ast),\\
            \hat \xi_h^\ast &\in\partial\TV_h(\hat u_h^\ast).
        \end{aligned}
    \right.
\end{equation}
Let us remark that it is straightforward to derive a version of \eqref{eq:kkt:gradient_h} in $Y_h\times Y_h\times U_h\times U_h\times\Psi_h$ instead of
$\R^{N_{Y_h}}\times\R^{N_{Y_h}}\times \R^{N_{U_h}}\times\R^{N_{U_h}}\times\R^{N_{\Psi_h}}$.
It can then be observed that this version is exactly \eqref{eq:optsys_pd} but with
$(\bar y,\bar p,\bar u,\bar q,\bar\xi)\in Y\times Y\times U\times U\times\Psi$ replaced by their finite-dimensional counterparts
$(y_h^\ast,p_h^\ast,u_h^\ast,q_h^\ast,\xi_h^\ast)\in Y_h\times Y_h\times U_h\times U_h\times\Psi_h$, and that \eqref{eq:kkt:gradient_h} is its equivalent reformulation in
$\R^{N_{Y_h}}\times\R^{N_{Y_h}}\times\R^{N_{U_h}}\times\R^{N_{U_h}}\times\R^{N_{\Psi_h}}$.
In particular, the two approaches of \emph{first discretize, then optimize} and \emph{first optimize, then discretize} coincide.

The next step is to characterize these subgradients componentwise. For the first subdifferential, we can simply use the sum and chain rules and find that
\begin{equation}
    [\hat q_h^\ast]_j\in d_j\partial\hat g([\hat u_h^\ast]_j),\qquad 1\leq j\leq N_{U_h},
\end{equation}
or equivalently
\begin{equation}
    [\hat u_h^\ast]_j \in\partial\hat g^*(d_j^{-1} [\hat q_h^\ast]_j), \qquad 1\leq j\leq N_{U_h},
\end{equation}
with $\partial\hat g^*$ given analogously to $\partial g^*$ as
\begin{equation}
    \partial \hat g^*(s) = \begin{cases}
        \{u_1\} & s \in \left(-\infty,\tfrac12(u_{1}+u_{2})\right),\\
        [u_i,u_{i+1}] & s = \tfrac12(u_{i}+u_{i+1}),\quad 1\leq i<m,\\
        \{u_i\} & s \in \left(\tfrac12(u_{i-1}+u_{i}),\tfrac12(u_{i}+u_{i+1})\right), \quad 1< i<m,\\
        \{u_m\} & s \in \left(\tfrac12(u_{m-1}+u_m),\infty\right),
    \end{cases}
\end{equation}
see also \cite[Sec.~2.1]{CK:2013}.
We will in the following replace the components $[\hat q_h^\ast]_i$ of $\hat q_h^\ast$ by their scaling $d_i^{-1} [\hat q_h^\ast]_i$; using the definition of the lumped mass matrix, this means we have to replace $\hat q_h^\ast$ in the third equation of \eqref{eq:kkt:gradient_h} by $M^\ell_h \hat q_h^\ast$.

For the discrete total variation, we use the sum rule and the chain rule
to deduce that there exists $\hat\psi_h^\ast\in\R^{N_{\Psi_h}}$ such that
\begin{equation}
    \hat\xi_h^\ast = D_h^T\hat\psi_h^\ast\qquad\quad\text{ and }\qquad\quad
    [\hat\psi_h^\ast]_T \in\partial(|\cdot|_2)([D_h\hat u_h^\ast]_T)\quad \text{for all }\,T\in\calT_h
\end{equation}
are satisfied.
As before, we rewrite the subdifferential inclusion equivalently as
\begin{equation}
    [D_h\hat u_h^\ast]_T \in\partial(|\cdot|_2^\ast) ([\hat\psi_h^\ast]_T)\quad \text{for all }\,T\in\calT_h.
\end{equation}
Using
\begin{equation}
    \hat h:\R^{d}\to\R,\qquad \hat h(v):= |v|_2,
\end{equation}
this reads
\begin{equation}
    [D_h\hat u_h^\ast]_T \in\partial\hat h^\ast ([\hat\psi_h^\ast]_T)\quad \text{for all }\,T\in\calT_h.
\end{equation}

\bigskip

To apply a Newton-type method, we replace the set-valued subdifferentials by their single-valued and Lipschitz-continuous Moreau--Yosida regularizations. Recall that the Moreau--Yosida regularization of $\partial F$ for any proper, convex and lower semi-continuous functional $F:X\to\Rbar:=\R\cup\{\infty\}$ acting on a Hilbert space $X$ is given by
\begin{equation}
    (\partial F)_\gamma(v) = \frac1\gamma\left(v-\prox_{\gamma F}(v)\right),
\end{equation}
where $\gamma>0$ and
\begin{equation}
    \prox_{\gamma F}(v) := \arg\min_{w\in
    X}\frac{1}{2\gamma}\norm{w-v}_{X}^2 + F(w)
    = \left(\Id + \gamma \partial F\right)^{-1}(v).
\end{equation}
For the regularized subdifferential $(\partial\hat g^*)_\gamma$, we have from \cite[Sec.~4.1]{CK:2015} that for $s\in\R$
\begin{equation}\label{eq:disc_reg_opt_q}
    (\partial \hat g^*)_\gamma(s) =
    \begin{cases}
        u_1 & s \in \left(-\infty,(\gamma+\tfrac12) u_{1}+\tfrac12 u_{2}\right),\\
        \frac1\gamma\left(s-\tfrac{u_i+u_{i+1}}{2}\right)& s \in \left[(\gamma+\tfrac12)u_i+\tfrac12u_{i+1},\tfrac12u_i +
        (\gamma+\tfrac12) u_{i+1}\right],\quad 1\leq i < m,\\
        u_i & s \in \left(\tfrac12 u_{i-1}+(\gamma+\tfrac12) u_{i},(\gamma+\tfrac12) u_{i}+u_{i+1}\right), \quad 1< i<m,\\
        u_m & s \in \left(\tfrac12 u_{m-1}+(\gamma+\tfrac12) u_m,\infty\right).
    \end{cases}
\end{equation}

For $\delta>0$, we denote the Moreau--Yosida regularization of $\partial\hat h^\ast$ by $(\partial \hat h^\ast)_\delta$.
To compute it, we recall that the Fenchel conjugate of a norm is the indicator function of the unit ball corresponding to the dual norm (which in this case is $|\cdot|_2$ itself). Furthermore, the proximal mapping $\prox_{\delta F}$ of an indicator function to a convex set is for every $\delta>0$ the metric projection onto this set. This shows that for all $v\in\R^d$ there holds
\begin{equation}
    (\partial \hat h^\ast)_\delta(v) =
    \frac1\delta\left(v - \proj_{\{|v|_2\leq 1\}}(v)\right) =
    \begin{cases}
        0 & |v|_2 \leq 1,\\
        \frac1\delta\left(v -\frac{v}{|v|_2}\right) & |v|_2 > 1.
    \end{cases}
\end{equation}

Combining the above, we obtain the regularized discrete optimality conditions
\begin{equation}\label{eq:disc_reg_opt_syshat}
    \left\{
        \begin{aligned}
            A_h(u_h^\ast+\umin) \hat y_h^\ast &= M_h \hat f_h,\\
            A_h(u_h^\ast+\umin) \hat p_h^\ast &= M_h(\hat z_h-\hat y_h^\ast),\\
            0&= \hat a_h(y_h^\ast,p_h^\ast) + \alpha M^\ell_h\hat q_h^\ast + \beta D_h^T \hat\psi_h^\ast,\\
            [\hat u_h^\ast]_j &= (\partial\hat g^*)_\gamma([\hat q_h^\ast]_j),\qquad\; 1\leq j\leq N_{U_h},\\
            [D_h\hat u_h^\ast]_T &= (\partial \hat h^\ast)_\delta([\hat\psi_h^\ast]_T),\qquad
            T\in \calT_h.
        \end{aligned}
    \right.
\end{equation}
Note that we have used the same notation $\hat y_h^\ast$, $\hat u_h^\ast$, etc., as for
solutions to the \emph{unregularized} discrete optimality conditions \eqref{eq:kkt:gradient_h} to avoid further complicating the notation. We point out that for the remainder of this work, this notation will always refer to solutions to \eqref{eq:disc_reg_opt_syshat}.

Finally, we remark that since $(\partial F^*)_\gamma = \nabla (F^*)_\gamma$ with $((F^*)_\gamma)^* = F + \frac\gamma2\norm{\cdot}^2_X$ holds for any proper, convex, and lower semi-continuous functional $F:X\to\Rbar$,
the regularized optimality system coincides with the necessary optimality conditions of
\begin{equation}\label{eq:problem_h_reg}
    \left\{\begin{aligned}
            \min_{\hat u_h\in\R^{N_{U_h}}} &\frac12\norm{y_h-z_h}_{L^2}^2 +
            \alpha\left(\hat G_h(\hat u_h)+\frac\gamma2\norm{\hat u_h}_{M^\ell_h}^2\right)
            +\beta\left(\TV_h(\hat u_h) + \frac\delta2 \norm{\hat u_h}_{2,h}^2\right)\\
            \text{s.t.}\quad &A_h(u_h+\umin)\hat y_h = M_h\hat f_h,
    \end{aligned}\right.
\end{equation}
where $\norm{\hat u_h}_{M^\ell_h}:=(\hat u_h^T M^\ell_h \hat u_h)^{1/2}$ and
$\norm{\hat u_h}_{2,h} := (\sum_{T\in\calT_h} | [D_h\hat u_h]_T |_2^2)^{1/2}$.
This can be interpreted as the mass-lumped approximation of an $H^1$ regularization of \eqref{eq:problem}.
Note, however, that the problem is still nonsmooth since $G_h$ and $\TV_h$ have not been modified; it has merely been made more strongly convex.

\subsection{A semismooth Newton-type method}\label{sec:semismooth}

To apply a semismooth Newton method to the regularized optimality conditions \eqref{eq:disc_reg_opt_syshat}, we reformulate them as a set of nonlinear implicit equations. Based on our numerical experience, it is preferable to consider the reduced system arising from \eqref{eq:disc_reg_opt_syshat} by eliminating the variables $(\hat u_h,\hat q_h)$ rather than solving the full system \eqref{eq:disc_reg_opt_syshat} in the variables $(\hat y_h,\hat p_h,\hat u_h, \hat q_h, \hat\psi_h)$.
In the following, we abbreviate $\hat\zeta_h:=(\hat y_h,\hat p_h,\hat\psi_h)\in\R^{N_{\hat\zeta_h}}$, where $N_{\hat\zeta_h}:=2 N_{Y_h}+N_{\Psi_h}$.

We begin the reformulation by noting that the third equation in \eqref{eq:disc_reg_opt_syshat} is equivalent to
\begin{equation}\label{eq:3rdofdiscreteKKTsystem}
    \hat q_h^\ast = - \frac1\alpha M^{-\ell}_h\left(B_h(y_h)\hat p_h + \beta D_h^T \hat\psi_h^\ast\right),
\end{equation}
where $M^{-\ell}_h$ denotes the inverse of $M^\ell_h$ and
$B_h(y_h)\in\R^{N_{U_h}\times N_{Y_h}}$ denotes the matrix induced by the bilinear form
\begin{equation}
    \scalprod{(\nabla y_h\cdot\nabla v_h),w_h} \qquad\text{for all }(w_h,v_h)\in U_h\times Y_h.
\end{equation}
Defining
\begin{equation}
    \hat q_h:\R^{N_{\hat\zeta_h}}\to\R^{N_{U_h}},
    \qquad \hat q_h(\hat\zeta_h):=-\frac1\alpha M^{-\ell}_h\left(B_h(y_h)\hat p_h + \beta D_h^T\hat\psi_h\right),
\end{equation}
\eqref{eq:3rdofdiscreteKKTsystem} becomes
\begin{equation}
    \hat q_h^\ast = \hat q_h(\hat\zeta_h^\ast).
\end{equation}
Inserting this into the fourth equation of \eqref{eq:disc_reg_opt_syshat} enables us to express $\hat u_h^\ast$ by
\begin{equation}
    \hat u_h^\ast = \hat u_h(\hat\zeta_h^\ast),
\end{equation}
where
\begin{equation}
    \hat u_h:\R^{N_{\hat\zeta_h}}\to\R^{N_{U_h}},
    \qquad
    \hat u_h(\hat\zeta_h):=
    \begin{pmatrix}
        (\partial \hat g^\ast)_\gamma ([\hat q_h(\hat \zeta_h)]_1)\\
        (\partial \hat g^\ast)_\gamma ([\hat q_h(\hat\zeta_h)]_2)\\
        \vdots\\
        (\partial \hat g^\ast)_\gamma ([\hat q_h(\hat\zeta_h)]_{N_{U_h}})
    \end{pmatrix}.
\end{equation}
We write $u_h(\hat\zeta_h)$ for the function $u_h\in U_h$ with coefficients $\hat u_h(\hat\zeta_h)$, i.e.,
$u_h(\hat\zeta_h):=\sum_{i=1}^{N_{U_h}} [\hat u_h(\hat\zeta_h)]_i\delta_i^{U_h}$.
Summarizing, \eqref{eq:disc_reg_opt_syshat} is equivalent to ${\mathcal F}_{\gamma,\delta}(\hat\zeta_h^\ast)=0$ for
\begin{equation}\label{calF}
    {\mathcal F}_{\gamma,\delta}:\R^{N_{\hat\zeta_h}}\to \R^{N_{\hat\zeta_h}},\qquad
    {\mathcal F}_{\gamma,\delta}(\hat\zeta_h) :=
    \begin{pmatrix}
        A_h(u_h(\hat\zeta_h)+\umin)\hat p_h + M_h (\hat y_h-\hat z_h)\\
        A_h(u_h(\hat\zeta_h)+\umin)\hat y_h - M_h \hat f_h\\
        {\cal H}(\hat\zeta_h)\\
    \end{pmatrix},
\end{equation}
where ${\cal H}:\R^{N_{\hat\zeta_h}}\to\R^{N_{\Psi_h}}$,
${\cal H}=({\cal H}_1^T,{\cal H}_2^T,\ldots,{\cal H}_{N_{\calT_h}}^T)^T$
with
\begin{equation}
    {\cal H}_j : \R^{N_{\hat{\zeta_h}}}\to \R^{d},
    \qquad
    {\cal H}_j(\hat\zeta_h):= [D_h \hat u_h(\hat\zeta_h)]_{T_j} - (\partial \hat h^\ast)_\delta([\hat\psi_h]_{T_j})\qquad \text{for }1\leq j\leq N_{\calT_h}.
\end{equation}
We recall that $\calT_h=\{T_j\}_{j=1}^{N_{\calT_h}}$ and point out
that $N_{\Psi_h} = N_{\calT_h} d$.

Since all components of ${\mathcal F}_{\gamma,\delta}$ are either
continuously differentiable or continuous and piecewise continuously differentiable
(PC$^1$) in each variable, ${\mathcal F}_{\gamma,\delta}$ is semismooth, see, e.g., \cite{Mifflin:1977,Kummer:1988,Kunisch:2008a,Ulbrich:2011}.
To obtain Newton derivatives for the nonsmooth terms, we use the fact that for PC$^1$ functions we can take as Newton derivative any selection of the derivatives of the essentially active pieces; see \cite[Sec.~2.5.3]{Ulbrich:2011}. In the following, we denote Newton derivatives by $D_N$. For the partial Newton derivative of, say, $\hat u_h(\cdot)$ with respect to the variable $\hat\psi_h$ evaluated at $\hat\zeta_h$, we write $D_{N_\psi}\hat u_h(\hat\zeta_h)$.
Since the mapping $\hat u_h(\cdot)$ is a composition of smooth mappings with $(\partial \hat g^*)_\gamma$, its Newton derivative is given by the chain rule in combination with our specific choice of
\begin{equation}
    D_N (\partial \hat g^*)_\gamma(s) =
    \begin{cases}
        \frac1\gamma & s \in\left[(\gamma+\tfrac12)u_i+\tfrac12u_{i+1},\tfrac12u_i +
        (\gamma+\tfrac12) u_{i+1}\right],\quad 1\leq i < m,\\
        0 & \text{else}.
    \end{cases}
\end{equation}
To determine $D_N{\cal H}$, it suffices to specify $D_N(\partial \hat h^\ast)_\delta$, where we make the choice
\begin{equation}
    D_N(\partial \hat h^\ast)_\delta(v) =
    \begin{cases}
        0 & |v|_2 \leq 1,\\
        \frac1\delta\left(\Id-\frac{1}{|v|_2}\Id+\frac1{|v|_2^3} v v^T\right) & |v|_2 > 1.
    \end{cases}
\end{equation}
Together, we obtain
\begin{equation}
    D_N {\cal F}_{\gamma,\delta}(\hat\zeta_h) =
    \begin{pmatrix}
        C_p E_y + M_h &
        C_p E_p + C_{y/p} &
        C_p E_\psi \\
        C_y E_y + C_{y/p} &
        C_y E_p &
        C_y E_\psi \\
        D_h E_y &
        D_h E_p &
        D_h E_\psi - E_{\psi\psi}
    \end{pmatrix}\in\R^{N_{\hat\zeta_h}\times N_{\hat\zeta_h}},
\end{equation}
where
\begin{align}
    C_p&:=B_h(p_h)^T, & C_y&:=B_h(y_h)^T, &
    C_{y/p}&:=A_h(u_h(\hat\zeta_h)+\umin),\\
    E_y&:=D_{N_y}\hat u_h(\hat\zeta_h), &
    E_p&:=D_{N_p}\hat u_h(\hat\zeta_h), &
    E_\psi&:=D_{N_\psi}\hat u_h(\hat\zeta_h),
\end{align}
and
\begin{equation}
    E_{\psi\psi}:=
    \begin{pmatrix}
        D_{N} (\partial \hat h^\ast)_\delta([\hat\psi_h]_{T_1}) & & & \\
                                                                & D_{N} (\partial \hat h^\ast)_\delta([\hat\psi_h]_{T_2}) & & \\
                                                                & & \ddots & \\
                                                                & & & D_{N} (\partial \hat h^\ast)_\delta([\hat\psi_h]_{T_{N_{\calT_h}}})
    \end{pmatrix}
    \in\R^{N_{\Psi_h}\times N_{\Psi_h}}.
\end{equation}

Note that the Newton matrix can become singular. For instance, if $|[\hat\psi_h]_T|_2 \leq 1$ for all $T\in\calT_h$, then $E_{\psi\psi}=0$. Hence, $(0,0,\hat w_h)^T\in\ker(D_N {\cal F}_{\gamma,\delta}(\hat\zeta_h))$ for every $\hat w_h\in\ker(E_\psi)$. Clearly, $\ker(E_\psi)$ is nontrivial since this is true for $\ker(D_h^T)$.
To cope with this singularity, we modify the (3,3) block of $D_N\calF_{\gamma,\delta}$ so that it reads $D_h E_\psi - E_{\psi\psi} - \mu_{\gamma,\delta}\hat M_h$, where $\hat M_h\in\R^{N_{\Psi_h}\times N_{\Psi_h}}$ denotes the diagonal mass matrix in $\Psi_h$, and $\mu_{\gamma,\delta}>0$ is a weight that depends on $\gamma$ and $\delta$; in our numerical experiments we observed $\mu_{\gamma,\delta}:=\delta^{-1}$ to work well. In the following, we assume that this choice is made unless explicitly indicated otherwise. We denote this modified matrix by $\widetilde{D_N\calF_{\gamma,\delta}}$.
For later reference we notice that given $(\gamma_j,\delta_j)\in\R_{>0}\times\R_{>0}$,
a semismooth Newton-type step $\tilde s^j\in\R^{N_{\hat\zeta_h}}$ at $\tilde\zeta^j\in\R^{N_{\hat\zeta_h}}$ is characterized by
\begin{equation}\label{eq:Newtonstep}
    \widetilde{D_N\calF_{\gamma_j,\delta_j}}(\tilde\zeta^j)\tilde s^j = -\calF_{\gamma_j,\delta_j}(\tilde\zeta^j).
\end{equation}
This step is combined with a backtracking line search based on the residual norm as well as a path-following scheme for $(\gamma_j,\delta_j)$.
The full procedure to compute an approximate solution to \eqref{eq:problem_hhat} is given in \cref{alg:TVMB},
where we have dropped the index $h$ for better readability.
We also write $\norm{\hat\zeta}_{L^2}:=\norm{\zeta_h}_{L^2(\Omega)^{d+2}}$ for $\hat\zeta\in\R^{N_{\hat\zeta}}$, where
$\hat\zeta$ are the coefficients of the function $\zeta_h\in Y_h\times Y_h\times\Psi_h$, and
$\norm{(\hat{\zeta},\hat u,\hat q)}_{L^2}:=\norm{(\zeta_h,u_h,q_h)}_{L^2(\Omega)^{d+4}}$,
where $(\hat{\zeta},\hat u,\hat q)$ are the coefficients of $(\zeta_h,u_h,q_h)\in Y_h\times Y_h\times\Psi_h\times U_h\times U_h$.

\begin{algorithm}[!t]
    \KwIn{
        \, $\hat\zeta^0\in\R^{N_{\hat\zeta}}$, \, $\gamma_0>0$, \, $\delta_0>0$, \, $\nu\in(0,1)$, \, $\text{TOL}_r>0$, \, $\text{TOL}_{\calF}>0$,\\
        \hspace{1.15cm} \, $\sigma_{\min}\in(0,1]$, \, $\sigma_{\text{nm}}\in(0,1]$
    }
    Set \, $k=0$ \, and \, $r_{-1}=\text{TOL}_r+1$\\
    \Repeat{$\bigl[\, r_{k-1}\leq\text{TOL}_r \bigr.$ \, \And \, $\bigl. r_{k-2}\leq \text{TOL}_r \,\bigr]$ }
    {
        Set \, $j=0$ \, and \, $\tilde\zeta^0=\hat\zeta^k$\\
        \While{\, $\norm{{\cal F}_{\gamma_j,\delta_j}(\tilde\zeta^j)}_{L^2} > \text{TOL}_{\calF}$ \,\label{alg:TVMB:4}}
        {
            Set \, $\tilde q^j=\hat q(\tilde\zeta^j)$ \, and \, $\tilde u^j=\hat u(\tilde\zeta^j)$\\
            Compute the Newton-type step \, $\tilde s^j$ \, at \, $\tilde\zeta^j$ \, by solving \eqref{eq:Newtonstep} and set \, $\sigma_j=1$\\
            \While{$\bigl[\, \sigma_j\geq\sigma_{\min} \bigr.$ \, \And \, $\bigl. \norm{\calF_{\gamma_j,\delta_j}(\tilde\zeta^j + \sigma_j \tilde s^j)}_{L^2}\geq\norm{\calF_{\gamma_j,\delta_j}(\tilde\zeta^j)}_{L^2} \,\bigr]$\label{alg:TVMB:7} }
            {
                Set \, $\sigma_j = \sigma_j/2$
            }
            \If{$\sigma_j<\sigma_{\min}$}{\, Set $\sigma_j = \sigma_{\text{nm}}$ \,}\label{alg:TVMB:12}
            Set \, $\tilde\zeta^{j+1} = \tilde\zeta^j + \sigma_j \tilde s^j$ \, and \, $j=j+1$
        }\label{alg:TVMB:14}
        Set \, $\hat\zeta^k_{\text{opt}} = \tilde\zeta^{j}$, \, $\hat u^{k}_{\text{opt}}=\hat u(\hat\zeta^k_{\text{opt}})$ \, and \, $\hat q^k_{\text{opt}}=\hat q(\hat\zeta^k_{\text{opt}})$\\
        Set \, $\gamma_{k+1} = \nu\gamma_k$ \, and \, $\delta_{k+1} = \nu\delta_k$\\
        \uIf{$k\geq 1$}
        {
            Set \, $r_k=\norm{(\hat\zeta^k_{\text{opt}},\hat u^k_{\text{opt}},\hat q^k_{\text{opt}})-(\hat\zeta^{k-1}_{\text{opt}},\hat u^{k-1}_{\text{opt}},\hat q^{k-1}_{\text{opt}})}_{L^2}$\\
            Set \, $\hat\zeta^{k+1}=\;\,(1+\nu)\hat\zeta^k_{\text{opt}}-\nu\hat\zeta^{k-1}_{\text{opt}}$\label{alg:TVMB:19}
        }
        \Else{Set \, $r_k=\text{TOL}_r+1$ \, and \, $\hat\zeta^{k+1}=\hat\zeta^k_{\text{opt}}$\label{alg:TVMB:21}}
        Set $k=k+1$
    }
    \KwOut{\, $\hat\zeta^{k-1}_{\text{opt}}\in\R^{N_{\hat\zeta}}$}
    \caption{Path-following method to solve \eqref{eq:problem_hhat}}
    \label{alg:TVMB}
\end{algorithm}

\Cref{alg:TVMB} is structured as follows. \Crefrange{alg:TVMB:4}{alg:TVMB:14} constitute an inner iteration; in this inner iteration, a Newton-type method with line search is employed for fixed $\gamma$ and $\delta$ to find a root of $\calF_{\gamma,\delta}$. The remaining lines form an outer iteration; in this outer iteration, $\gamma$ and $\delta$ are updated and the starting point for the next inner iteration is computed in \cref{alg:TVMB:19} or \cref{alg:TVMB:21}, respectively. Moreover, the $L^2$ difference of subsequent outer iterates is stored in $r_k$ and used in the termination criterion.

Let us comment on some important features of \cref{alg:TVMB}. We start by pointing out that the line search in \crefrange{alg:TVMB:7}{alg:TVMB:12} of \cref{alg:TVMB} is nonmonotone. That is, if backtracking does not yield a $\sigma_j\in[\sigma_{\min},1]$ with $\norm{\calF_{\gamma_j,\delta_j}(\tilde\zeta^j + \sigma_j\tilde s^j)}_{L^2}<\norm{\calF_{\gamma_j,\delta_j}(\tilde\zeta^j)}_{L^2}$,
then the step length $\sigma_j=\sigma_{\text{nm}}$ is used regardless whether it satisfies
$\norm{\calF_{\gamma_j,\delta_j}(\tilde\zeta^j + \sigma_{\text{nm}} \tilde s^j)}_{L^2}<\norm{\calF_{\gamma_j,\delta_j}(\tilde\zeta^j)}_{L^2}$ or not.

Next we remark that the computation of $\hat\zeta_{k+1}$ in \cref{alg:TVMB:19} is a predictor step: From the previous roots $\hat\zeta_{\text{opt}}^k$ and $\hat\zeta_{\text{opt}}^{k-1}$, a prediction
$\hat\zeta^{k+1}$ of $\hat\zeta_{\text{opt}}^{k+1}$ is computed and used as the starting point for the next inner iteration (whose aim it is to find $\hat\zeta_{\text{opt}}^{k+1}$). For $k\geq 1$, this prediction is taken to be the componentwise linear extrapolation
\begin{equation}
    \hat\zeta^{k+1}:= \hat\zeta_{\text{opt}}^{k} + \frac{\gamma_k-\gamma_{k+1}}{\gamma_{k-1}-\gamma_k} (\hat\zeta_{\text{opt}}^{k}-\hat\zeta_{\text{opt}}^{k-1}) = 
    (1+\nu)\hat\zeta^k_{\text{opt}}-\nu\hat\zeta^{k-1}_{\text{opt}},  
\end{equation}
where we have used that $\gamma_{k+1}=\nu\gamma_k=\nu^2\gamma_{k-1}$.
Note that due to the coupling $\delta_k/\gamma_k=\delta_0/\gamma_0$ for all $k$, we obtain the same extrapolation step if $\gamma$ is replaced by $\delta$. We thus perform a combined prediction for the continuation in $\gamma$ as well as $\delta$. 
For $k=0$, no predictor step is used as
$\hat\zeta^{k-1}_{\text{opt}}=\hat\zeta^{-1}_{\text{opt}}$ is not available; instead we set $\hat\zeta^1 = \hat\zeta^{0}_{\text{opt}}$ in this case.

Finally, we embed \cref{alg:TVMB} within a further continuation strategy for $\nu$: If a Newton iteration for a given pair $(\gamma_k,\delta_k)$ does not terminate successfully, we increase $\nu$ and restart \cref{alg:TVMB} from the last successful solution; this outer continuation is terminated if $\nu\approx 1$.

\bigskip

We conclude this section with several practical remarks concerning \cref{alg:TVMB}.
First, we stress that while its numerical costs are negligible, the	
predictor step significantly increased the convergence speed in our numerical experiments.
Also, due to the path-following strategy, it is not necessary to choose the initial guess $\hat\zeta^0$ in a specific way. In fact, our numerical experiments indicate that arbitrary starting points can be used.
In particular, the choice $\hat\zeta^0:=0$ was always sufficient to achieve convergence.

Furthermore, we found in our numerical experiments that for larger values of $\gamma$ and $\delta$ (e.g., $\gamma,\delta>1$), the convergence of \cref{alg:TVMB} can be accelerated if
$\mu_{\gamma,\delta}=\delta$ is used and $\calF_{\gamma,\delta}$ is modified such that its Newton derivative equals $\widetilde{D_N\calF_{\gamma,\delta}}$.
For small values of $\gamma$ and $\delta$, however, this strategy did not work and we had to choose $\calF_{\gamma,\delta}$ as given in \eqref{calF} and $\mu_{\gamma,\delta}=\delta^{-1}$. Note that for the choice $\mu_{\gamma,\delta}=\delta^{-1}$ it is not sensible to modify $\calF_{\gamma,\delta}$ in such a way that its Newton derivative equals $\widetilde{D_N\calF_{\gamma,\delta}}$. In fact, we can show that if $\calF_{\gamma,\delta}$ is modified in this way, then the sequence $((\hat\zeta_{\text{opt}}^k,\hat u(\hat\zeta^k_{\text{opt}}),\hat q(\hat\zeta^k_{\text{opt}})))_k$
can only converge to a solution to \eqref{eq:kkt:gradient_h} with $\beta=0$, i.e., to a solution to the optimality conditions of the ``pure multi-bang problem''.

\section{Numerical examples}\label{sec:examples}

We illustrate the structure of optimal controls for \eqref{eq:problem_hhat} using two model problems. In particular, the goal is to show the difference between optimal controls of \eqref{eq:problem_hhat} for $\beta>0$ and for $\beta=0$, i.e., between solutions to a TV-regularized multi-bang problem and those to a ``pure multi-bang'' problem. We remark that $\beta>0$ is required in the infinite dimensional case but can be arbitrarily small, while taking $\beta=0$ is justified in the finite-dimensional setting only. More examples for the pure multi-bang approach can be found in \cite{CK:2013,CK:2015}.

In all examples, we take $\Omega=(-1,1)^2\subset\R^2$ and employ a uniform triangulation $\calT_h$ consisting of $8192$ elements, i.e., $N_{U_h}=64\cdot 64$. We use
$\umin=1.5$ and the algorithmic parameters
$\hat\zeta^0 = 0$, $\gamma_0 = 10^5$, $\delta_0 =10^3$, $\nu=0.8$, $\nu_{\max}=0.9999$, $\text{TOL}_r = 10^{-3}(\umax-\umin)$, $\text{TOL}_{\calF} = 10^{-5}$, as well as
$\sigma_{\min}=10^{-6}$ and $\sigma_\text{nm}=10^{-2}$.
The remaining data and parameters are chosen individually for each example.

We implemented \cref{alg:TVMB} in Python using DOLFIN \cite{LoggWells2010a,LoggWellsEtAl2012a}, which is part of the open-source computing platform FEniCS \cite{AlnaesBlechta2015a,LoggMardalEtAl2012a}. The linear system \eqref{eq:Newtonstep} arising from the Newton-type step is solved using the sparse direct solver \verb!spsolve! from SciPy.

\subsection{Example 1: topology optimization}

The first example is motivated by the possible application to topology optimization. The general idea is that we have a design $\tilde u\in U$ making use of two materials characterized by their densities $\umin+\tilde u_1=1.5$ and $\umin+\tilde u_2=2.5$; we call this a \emph{binary design}. Imagine that it has become possible to use also materials that have intermediate densities, e.g., in total five materials with
densities $\umin+u_j=1.5+0.25(j-1)$, $1\leq j\leq 5$. The question is now whether it is possible to realize a similar state as arising from the (presumably optimal) binary design using the (presumably cheaper) intermediate materials. 

\begin{figure}[t]
    \centering
    \begin{subfigure}[t]{0.495\linewidth}
        \centering
        \includegraphics[width=\textwidth]{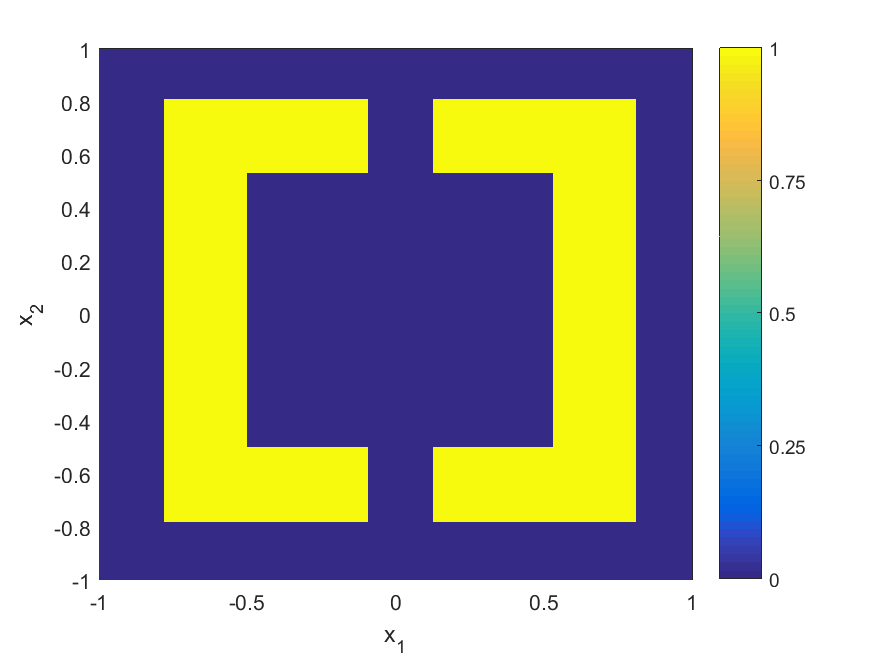}
        \caption{binary design $\tilde u_h-\umin$}\label{fig:ex1:binary}
    \end{subfigure}
    \hfill
    \begin{subfigure}[t]{0.495\linewidth}
        \centering
        \includegraphics[width=\textwidth]{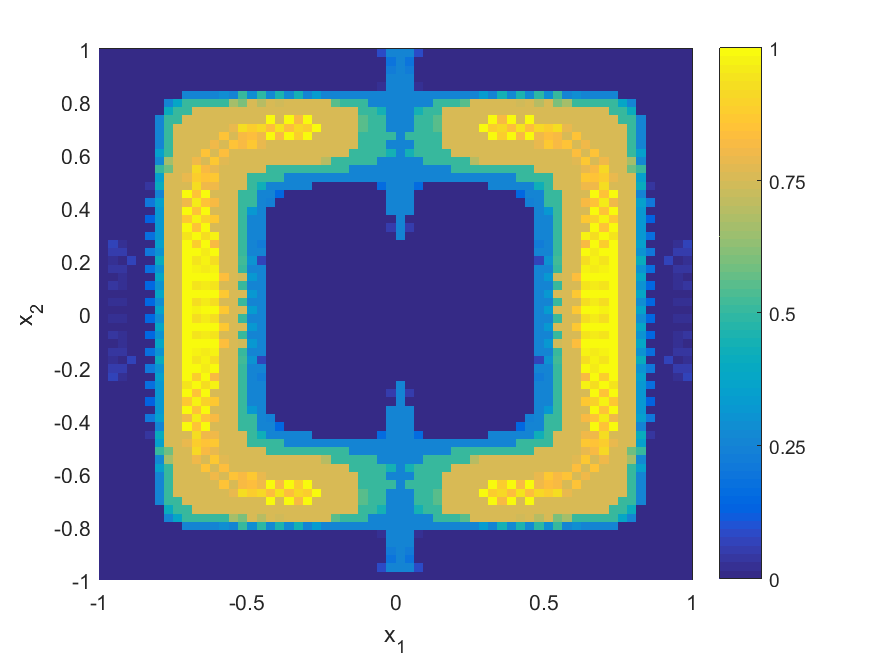}
        \caption{pure multi-bang design $\bar u_h-\umin$}\label{fig:ex1:multi-bang}
    \end{subfigure}

    \begin{subfigure}[t]{0.495\linewidth}
        \centering
        \includegraphics[width=\textwidth]{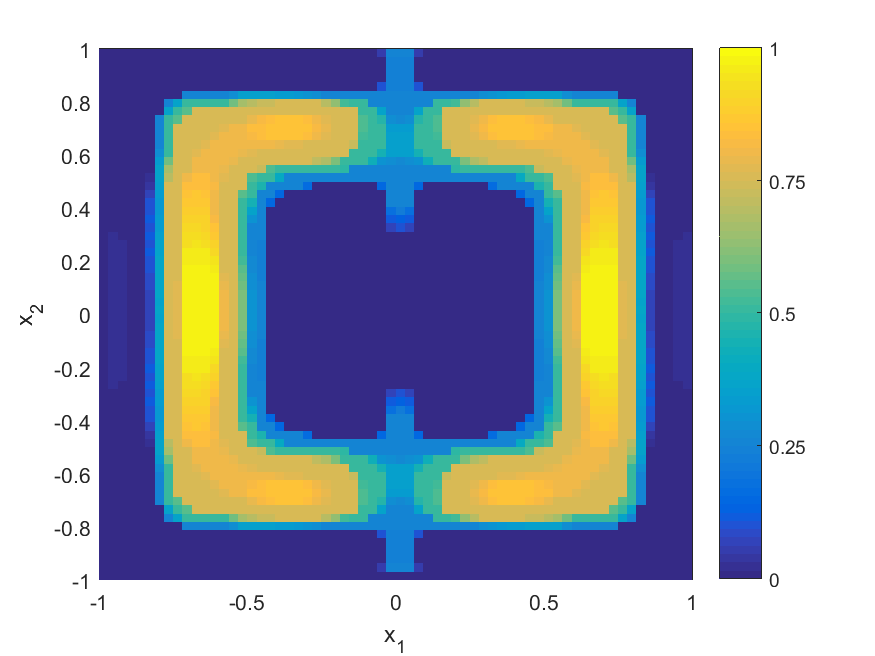}
        \caption{TV--multi-bang design $u_h^\ast-\umin$, $\beta=10^{-6}$}\label{fig:ex1:tv}
    \end{subfigure}
    \hfill
    \begin{subfigure}[t]{0.495\linewidth}
        \centering
        \includegraphics[width=\textwidth]{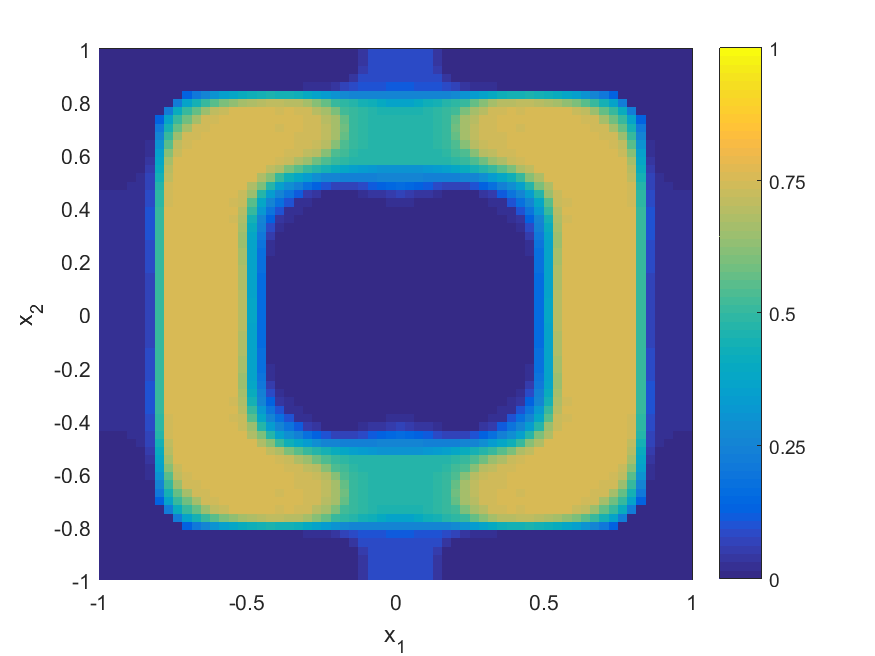}
        \caption{TV--multi-bang design $u_h^\ast-\umin$, $\beta=5\cdot 10^{-5}$}\label{fig:ex1:tvlargebeta}
    \end{subfigure}
    \caption{Comparison of binary, pure multi-bang, and total variation designs for Example 1}
    \label{fig:ex1}
\end{figure}

Following this motivation, we start from the binary design
\begin{equation}
    \tilde u(x):=\begin{cases}
        1.5 & x\in\omega_1,\\
        2.5 & x\in\omega_2,
    \end{cases}
\end{equation}
where
\begin{equation}
    \omega_2 :=\Bigl\{
    x\in\Omega: 0.1<|x_1|<0.8 \quad \text{and}\quad |x_2|<0.8\quad\text{and}\quad \bigl[|x_1|>0.5 \,\text{ or }\, |x_2|>0.5\bigr]\Bigr\}
\end{equation}
and $\omega_1:=\Omega\setminus\omega_2$.
Denoting by $\tilde u_h\in U_h$ the finite element function that interpolates $\tilde u$ in all vertices of $\calT_h$, we compute the target $z_h\in Y_h$ as the state corresponding to $\tilde u_h$ and $f_h\equiv 10$, i.e., as the solution to $-\div(\tilde u_h \nabla z_h) = f_h$ in $\Omega$; see \cref{fig:ex1:binary}.
We then compute a solution to \eqref{eq:problem_hhat} using the five desired coefficient values $u_j =0.25(j-1)$, $1\leq j\leq 5$, together with the parameters
$\alpha=10^{-3}$ and $\beta\in\{0,10^{-6},5\cdot10^{-5}\}$; see \crefrange{fig:ex1:multi-bang}{fig:ex1:tvlargebeta} (with $\gamma_{\text{final}}\approx 1\cdot 10^{-4}$, $\gamma_{\text{final}}\approx 1.8\cdot 10^{-4}$, and
$\gamma_{\text{final}}\approx 6.4\cdot 10^{-2}$, respectively).

Comparing the pure multi-bang design $\bar u_h$ in \cref{fig:ex1:multi-bang} with the TV--multibang designs in \cref{fig:ex1:tv}--\subref{fig:ex1:tvlargebeta}, we clearly observe the well-known effect of TV regularization favoring level sets with smaller perimeter: While most jumps and the promotion of the desired parameter values are retained from the pure multi-bang design, the high-frequency ``oscillations'' between the level sets of $\bar u_h = 1.0$ and $\bar u_h = 0.75$ are removed. Similarly, the spurious ``droplets'' near $x=(-1,0)$ and $x=(1,0)$ are suppressed. (Here we recall that the multi-bang penalty acts purely pointwise and does not promote any spatial regularity.) The effect of the total variation penalty is also visible in \cref{fig:ex1:tvlargebeta}, where the perimeters of the level sets for $u_h^\ast=0.5$ and $u_h^\ast = 0.75$ have both been reduced, respectively, by closing the ``slit'' at $x_1=0$ and by removing the highest-valued material.
We point out that the simpler structure of the TV-regularized control may in itself be preferable in certain applications.
(We also remark that if the admissible control values are restricted to $(u_1,u_2)=(\tilde u_1,\tilde u_2)=(0,1)$ and $\alpha,\beta$ are chosen sufficiently small, the binary reference design is essentially recovered.)

\subsection{Example 2: parameter identification}

The second example is motivated by a parameter identification related to electrical impedance tomography. Here, the goal is to reconstruct the spatially varying conductivity (which is a tissue-specific material parameter) from noisy observations of the electric field arising from external charges. It should be noted that in medical impedance tomography, external currents and observations are both taken on the boundary or a part thereof; for the sake of simplicity, however, we consider distributed charge density and observation.

\begin{figure}[t]
    \centering
    \begin{subfigure}[t]{0.495\linewidth}
        \centering
        \includegraphics[width=\textwidth]{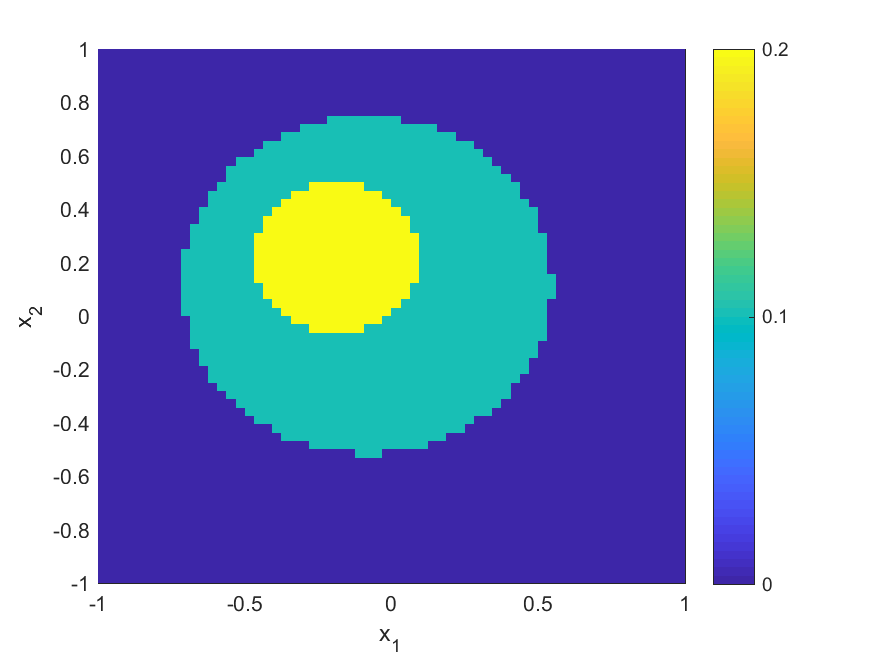}
        \caption{true parameter $\tilde u_h-\umin$}\label{fig:ex2:binary}
    \end{subfigure}
    \hfill
    \begin{subfigure}[t]{0.495\linewidth}
        \centering
        \includegraphics[width=\textwidth]{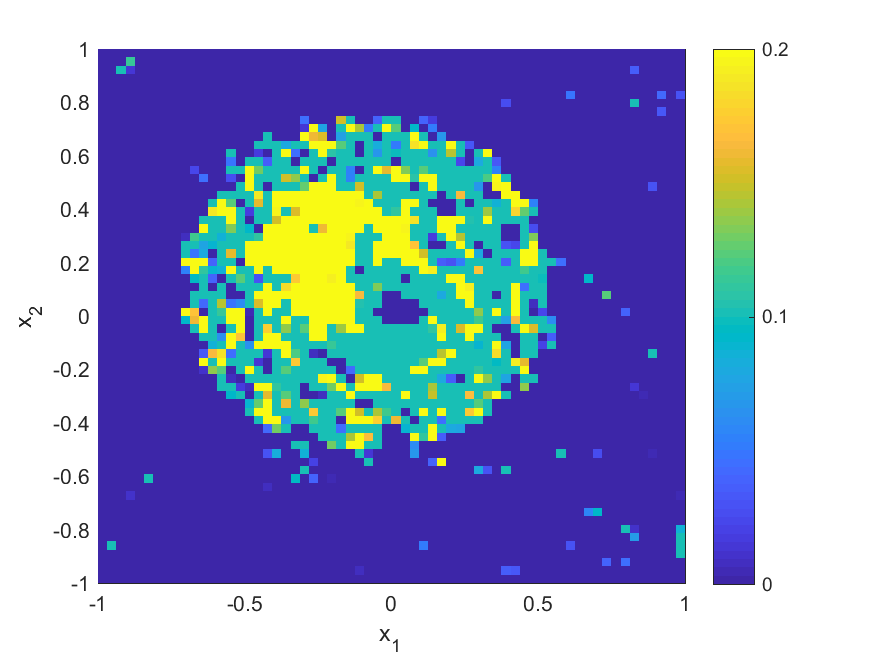}
        \caption{multi-bang reconstruction $\bar u_h-\umin$}\label{fig:ex2:multi-bang}
    \end{subfigure}

    \begin{subfigure}[t]{0.495\linewidth}
        \centering
        \includegraphics[width=\textwidth]{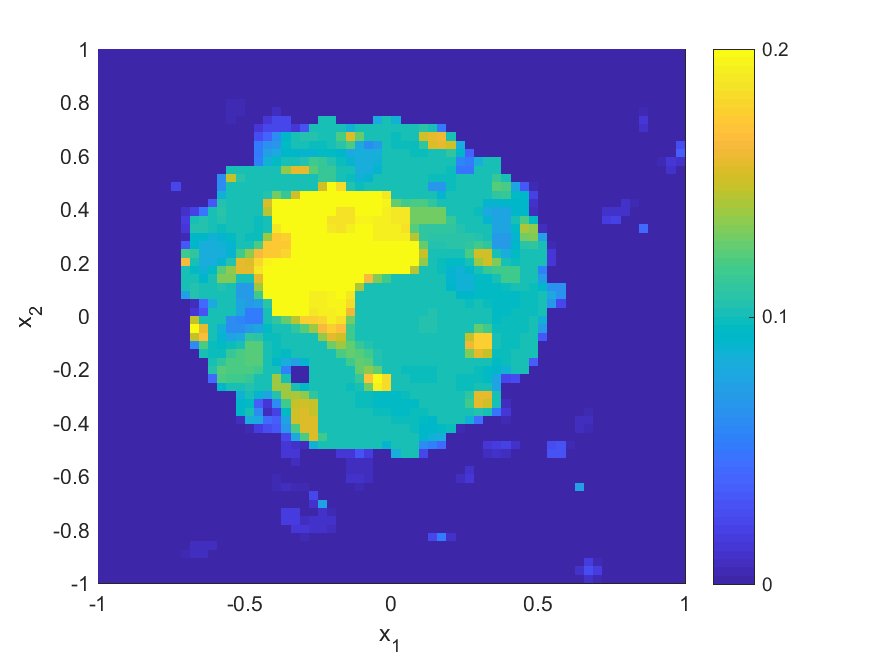}
        \caption{TV--multi-bang recon. $u_h^\ast-\umin$, $\beta=10^{-6}$}\label{fig:ex2:tv_small}
    \end{subfigure}
    \hfill
    \begin{subfigure}[t]{0.495\linewidth}
        \centering
        \includegraphics[width=\textwidth]{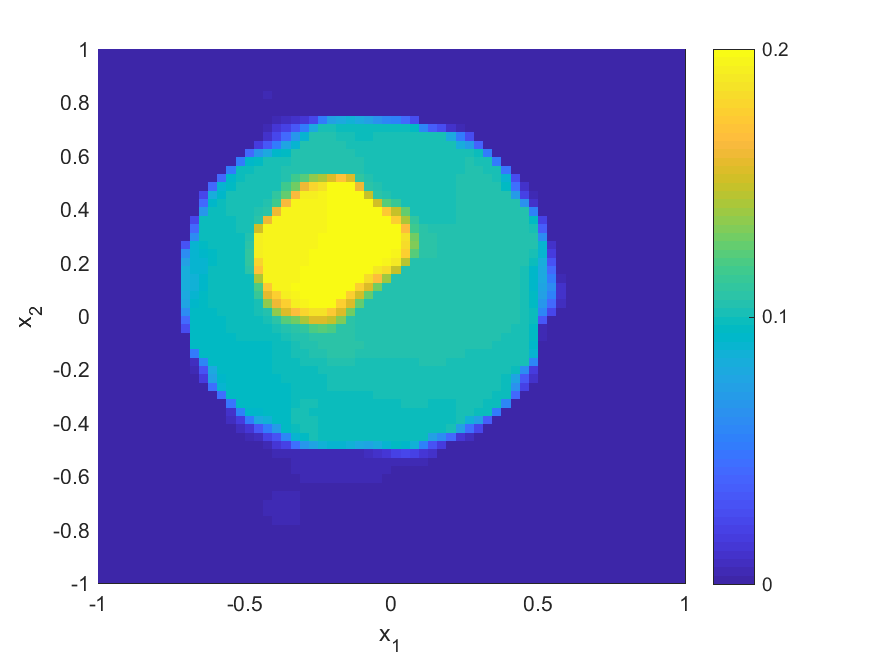}
        \caption{TV--multi-bang recon. $u_h^\ast-\umin$, $\beta=10^{-5}$}\label{fig:ex2:tv}
    \end{subfigure}
    \caption{Comparison of true parameter, pure multi-bang, and total variation-regularized reconstructions for Example 2}
    \label{fig:ex2}
\end{figure}

We choose as true parameter
\begin{equation}
    \tilde u(x):=\begin{cases}
        1.5 & x\in\omega_1,\\
        1.6 & x\in\omega_2,\\
        1.7 & x\in\omega_3,
    \end{cases}
\end{equation}
where
\begin{equation}
    \omega_1 := \Bigl\{x\in\Omega: (x_1+0.1)^2+(x_2-0.1)^2\geq 0.4\Bigr\},\qquad
    \omega_3 := \Bigl\{x\in\Omega: (x_1+0.2)^2+(x_2-0.2)^2<0.08\Bigr\},
\end{equation}
and $\omega_2:=\Omega\setminus(\omega_1\cup\omega_3)$ model background, tumor, and healthy tissue, respectively.
Again, $\tilde u_h\in U_h$ denotes the finite element function interpolating $\tilde u$ in all vertices of $\calT_h$; see \cref{fig:ex2:binary}.
For the target, we first compute a noise-free state $\tilde z_h\in Y_h$ solving $-\div(\tilde u_h \nabla \tilde z_h) = f_h$ in $\Omega$, where $f_h\equiv 25$. We now add noise to $\tilde z_h$ to obtain $z_h$; we use $z_h:=\tilde z_h+n_l\rho_h\max_{x\in\Omega}(|\tilde z_h(x)|) $, where $n_l:=10^{-3}$ and $\rho_h\in Y_h$ is a finite element function whose coefficients $\hat\rho_h\in\R^{N_{Y_h}}$ are sampled from a normal distribution with mean zero and standard deviation one.
Corresponding to the assumption that strong a priori knowledge is available, we choose
the desired coefficient values $u_1=0$, $u_2=0.1$ and $u_3=0.2$, together with the parameters
$\alpha=5\cdot 10^{-4}$ and $\beta\in\{0,10^{-5},10^{-6}\}$; see \crefrange{fig:ex2:multi-bang}{fig:ex2:tv} (with $\gamma_{\text{final}}\approx 5.8\cdot 10^{-6}$, $\gamma_{\text{final}}\approx 2.9\cdot 10^{-3}$, and
$\gamma_{\text{final}}\approx 6.6\cdot 10^{-3}$, respectively).

From \cref{fig:ex2:multi-bang}, it is obvious that the pure multi-bang regularization fails for this challenging problem since the multi-bang penalty entails no spatial regularization. Specifically, noise remains in the homogeneous background, and many points in the healthy tissue region are misclassified as either tumor or background; the latter in particular in a large region near $x=(0,0)$ where $\nabla \bar y_h\approx 0$ (compare \eqref{eq:kkt}). The reconstruction is improved by adding the total variation regularization: with $\beta=10^{-6}$, the ``hole'' near $x=(0,0)$ is gone, and the misclassified points are reduced; see \cref{fig:ex2:tv_small}. Increasing the total variation regularization parameter to $\beta=10^{-5}$ (\cref{fig:ex2:tv}) again significantly improves the reconstruction by removing the small spurious inclusions while preserving the contrast and shape of the healthy tissue and tumor regions; merely the volume of the latter is slightly reduced. This indicates that regularization as understood in the context of inverse problems is predominantly provided by the total variation penalty, while the multi-bang penalty is responsible for maintaining the desired contrast of the reconstruction. Hence, it suffices to investigate noise level-dependent parameter choice rules for $\beta$ while keeping $\alpha$ fixed, rather than having to consider -- much more challenging -- choice rules for multiple parameters.

\section{Conclusion}

Total variation regularization of topology optimization and parameter identification problems is challenging both analytically and numerically but is required in order to obtain existence of a solution without introducing additional smoothing. Furthermore, a pointwise multi-bang penalty can be used to promote optimal coefficients with desired (material) values. A reparametrization of the coefficient to be optimized allows proving existence as well as obtaining pointwise optimality conditions.
The numerical solution is based on a finite element discretization and Moreau--Yosida regularization of reduced optimality conditions together with a semismooth Newton-type method combined with a predictive path-following strategy.
Numerical examples indicate that in comparison to a pure multi-bang approach, the additional total variation regularization yields controls whose structure is much more regular. 

\appendix

\section{Strongly Lipschitz domains are Gröger regular}\label{sec:appendix}

In this appendix, we address the relation between two different definitions of Lipschitz domains and the concept of Gröger regularity which are used in the literature. The first definition, sometimes referred to as a \emph{strongly Lipschitz domain}, requires that, roughly speaking, the boundary can be represented locally as the graph of a Lipschitz function. A precise statement is the following from \cite[A\,8.2]{Alt:EnglishVersion}.
\begin{definition}[Strongly Lipschitz domain]\label{def:Lipboundary_Alt}
    Let $\Omega\subset\R^d$ be open and bounded. We say that $\Omega$ has a \emph{Lipschitz boundary} if there exists $l\in\N$ 
    such that $\partial\Omega$ can be covered by open sets $U^1,U^2,\ldots,U^l$
    and for $j=1,\ldots,l$ there exist a Euclidean coordinate system $e_1^j,e_2^j,\ldots,e_d^j\in\R^d$, a reference point 
    $y^j\in\R^{d-1}$, numbers $r^j>0$ and $h^j>0$, and a Lipschitz continuous function $\eta^j:\R^{d-1}\rightarrow\R$ that satisfy the following properties:
    \begin{enumerate}[label=(\roman*)]
        \item $U^j = \left\{x\in\R^d: |x_{-,d}^j-y^j|_2<r^j \; \text{ and } \; |x_d^j-\eta^j(x_{-,d}^j)|_2<h^j\right\}$;
        \item for all $x\in U^j$, if $x_d^j=\eta^j(x_{-,d}^j)$ then $x\in\partial\Omega$;
        \item for all $x\in U^j$, if $0<x_d^j-\eta^j(x_{-,d}^j)<h^j$ then $x\in\Omega$;
        \item for all $x\in U^j$, if $0>x_d^j-\eta^j(x_{-,d}^j)>-h^j$ then $x\not\in\Omega$.
    \end{enumerate}	
    Here, we have denoted $x_{-,d}^j = (x_1^j,\ldots,x_{d-1}^j)^T\in\R^{d-1}$ for $x=x^j=(x_1^j,\ldots,x_d^j)^T\in\R^d$, and the coordinates of $x^j$ are given in the local Euclidean coordinate system $e_1^j,e_2^j,\ldots,e_d^j$ in $\R^d$, i.e., $x^j = \sum_{i=1}^d x_i^j e_i^j$.	

    A bounded domain with Lipschitz boundary is called a \emph{strongly Lipschitz domain}.
\end{definition}
Strongly Lipschitz domains are extension domains, which is required to obtain embeddings for Sobolev and BV functions into $L^p$ spaces, and this definition is therefore used in \cite{Ambrosio,BrediesHoller}.

The second definition, sometimes referred to as a \emph{weakly Lipschitz domain}, requires, roughly speaking, that the boundary can be locally flattened by a bi-Lipschitz transformation; a precise definition can be found in, e.g., \cite[Sec.~6]{Brewster}.
For our purposes, however, the following related concept from \cite[Def.~2]{Groeger} is more important.
\begin{definition}[Gröger regularity]\label{def:Lipboundary_Groeger}
    A set $G\subset\R^d$ is called \emph{regular (in the sense of Gröger)} if $G$ is bounded and if for every $y\in\partial G$ there exist subsets $U$ and $\tilde U$ of $\R^d$ and a Lipschitz continuous bijection $\Phi: U\rightarrow\tilde U$ with Lipschitz continuous inverse $\Phi^{-1}$ such that $U$ is an open neighborhood of $y$ in $\R^d$ and that $\Phi(U\cap G)$ is one of the sets
    \begin{align}
        E_1&:=\left\{x\in\R^d: |x|<1,\,x_d<0\right\},\\
        E_2&:=\left\{x\in\R^d: |x|<1,\,x_d\leq 0\right\},\\
        E_3&:=\bigl\{x\in E_2: x_d<0 \; \text{ or } \; x_1>0\bigr\},
    \end{align}
    where $x=(x_1,x_2,\ldots,x_d)^T$.
\end{definition}
The main result of \cite{Groeger} is that a second order elliptic mixed boundary value problem on a bounded domain $\Omega$ admits higher regularity of the solution if $G=\Omega\cup \Gamma_N$ is regular, where $\Gamma_N\subset \partial\Omega$ denotes the Neumann boundary. For $\Gamma=\emptyset$ (i.e., pure Dirichlet conditions, where $E_2$ and $E_3$ are not needed), \cref{def:Lipboundary_Groeger} reduces to that of $\Omega$ being a weakly Lipschitz domain. Furthermore, \cite[Sec.~5]{HDMRS:2009} shows also for mixed boundary conditions (under some assumptions on $\Gamma_N$) that if $\Omega\cup \Gamma_N$ is regular then $\Omega$ is a weakly Lipschitz domain and, for $d\in\{2,3\}$, vice versa.

In our analysis, we require the domain $\Omega$ to satisfy both \cref{def:Lipboundary_Alt} and \cref{def:Lipboundary_Groeger} since we use results from \cite{Ambrosio,BrediesHoller} as well as from \cite{Groeger}.
However, the notions of strongly and weakly Lipschitz domains are not equivalent; examples of weakly but not strongly Lipschitz domains can be found in, e.g., \cite[Sec.~6]{Brewster}.
Although it is commonly accepted that strongly Lipschitz domains are regular (or, equivalently for domains, that they are weakly Lipschitz), despite our best efforts we could not find a proof of this fact in the literature. For the sake of completeness, we therefore provide one here.
\begin{lemma}\label{prop:regularity:LipschitzimpliesGroegerregularity}
    If a domain $\Omega\subset\R^d$ satisfies \cref{def:Lipboundary_Alt}, then it also satisfies \cref{def:Lipboundary_Groeger}.
\end{lemma}
\begin{proof}
    Let $\Omega\subset\R^d$, $d\in\N$, denote the set in question and let $\hat x\in\partial\Omega$. Due to \cref{def:Lipboundary_Alt} there exist an open neighborhood $\hat V$ of $\hat x$ and a Lipschitz continuous function $\eta:\R^{d-1}\rightarrow\R$ such that
    $\Omega\cap\hat V = \{x\in\hat V: \eta(x_1,\ldots,x_{d-1})<x_d\}$ and $\eta(\hat x_1,\ldots,\hat x_{d-1})=\hat x_d$.
    Defining
    \begin{equation}
        \hat\Lambda:\R^d\rightarrow\R^d,\qquad \hat\Lambda(x):=(x_1,\ldots,x_{d-1},\eta(x_1,\ldots,x_{d-1})-x_d)^T,
    \end{equation}
    we observe that $\Omega\cap\hat V = \{x\in\hat V: \hat\Lambda_d(x)<0\}$.
    Clearly, $\hat\Lambda$ is Lipschitz. Moreover, since $\hat\Lambda(\hat\Lambda(x)) = x$ for all $x\in\R^d$,
    we infer that $\hat\Lambda$ and its inverse mapping $\hat\Xi:=\hat\Lambda^{-1}=\hat\Lambda$ are bijective. (In the following, we nevertheless distinguish between $\hat\Lambda$ and its inverse for the sake of transparency.)
    Since $\hat\Lambda$ is Lipschitz continuous, $\hat\Xi$ maps open sets to open sets.
    Defining $\hat y:=\hat\Lambda(\hat x)$ we note that $\hat\Xi$ maps $B_\delta(\hat y)$ for every $\delta>0$ bijectively to $\hat\Xi(B_\delta(\hat y))$, which is an open neighborhood of $\hat x$.
    In particular, there is $\delta>0$ such that $\hat\Xi$ maps $B_\delta(\hat y)$ bijectively to $V:=\hat\Xi(B_\delta(\hat y))$ with $V\subset\hat V$.
    Consequently, $\Xi(y):=\hat\Xi(\hat y+\delta y)$ maps $B_1(0)$ bijectively to $V$, is Lipschitz continuous, and has the Lipschitz continuous inverse $\Lambda(x):=(\hat\Lambda(x)-\hat y)/\delta$.
    It follows that
    \begin{equation}
        \Omega\cap V=\{x\in V: \hat\Lambda_d(x)<0\}=\{\Xi(y)\in\R^d: y\in B_1(0),\, \hat\Lambda_d(\Xi(y))<0\}.
    \end{equation}
    This implies that
    \begin{equation}
        \Lambda(\Omega\cap V)=\{y\in\R^d: y\in B_1(0),\,\hat y_d +\delta y_d<0\}=\{y\in\R^d: y\in B_1(0),\, y_d<0\},
    \end{equation}
    where we have used $\hat y_d=0$.
    Summarizing, we have established that for $\hat x\in\partial\Omega$, there is an open neighborhood $V$ of $\hat x$ and a Lipschitz continuous bijection $\Lambda:V\rightarrow B_1(0)$ with Lipschitz continuous inverse such that
    $\Lambda(\Omega\cap V)=\{y\in B_1(0): y_d<0\}$. That is, $G:=\Omega$ satisfies \cref{def:Lipboundary_Groeger}.
\end{proof}

\section*{Acknowledgments}

Support by the German Science Fund (DFG) under grant CL 487/1-1 for C.C. and by the ERC advanced grant 668998 (OCLOC) under the EU's H2020 research program for K.K. are gratefully acknowledged.

\printbibliography

\end{document}